\documentclass[a4paper,reqno]{amsart}
\usepackage{amsthm,amsfonts,amsmath,enumerate,amssymb,stmaryrd} 
\usepackage{graphics}                 
\usepackage{color}                    
\usepackage{cancel}
\usepackage{todonotes}
\usepackage{johncmd}

\usepackage{multirow}

\usepackage{xr-hyper}
\usepackage[colorlinks=true]{hyperref}
\usepackage{verbatim}
\usepackage{pdfsync}

\usepackage[textheight=630pt,
  textwidth=468pt,
  centering]{geometry}

\newcommand{\myeq}[2]{\stackrel{\makebox[0pt]{\mbox{\normalfont\small #2}}}{#1}}

\newcommand{\cs}[1]{{#1}}

\newcommand{\changes}[1]{{#1}}

\renewcommand{\Re}{{\mathrm Re}}
\newcommand{\dfrak}{\Delta}

\newcommand{\deltau}[1]{\dfrak #1}
\newcommand{\ucs}[1]{\widehat{{#1}}}
\newcommand{\ii}{\boldsymbol{i}}
\newcommand{\bsnul}{\boldsymbol{0}}
\newcommand{\IR}{\mathbb R} 
\newcommand{\IC}{\mathbb C} 
\newcommand{\IN}{\mathbb N} 
\newcommand{\bcX}{\mathcal{X}}
\newcommand{\bcY}{{\mathcal{Y}}}
\newcommand{\bas}{\begin{align*}}
\newcommand{\eas}{\end{align*}}
\newcommand{\ba}{\begin{align}}
\newcommand{\ea}{\end{align}}
\newcommand{\bes}{\begin{equation*}}
\newcommand{\ees}{\end{equation*}}
\newcommand{\be}{\begin{equation}}
\newcommand{\ee}{\end{equation}}

\newcommand{\Tscheb}{Chebyshev\xspace}

\newcommand{\cL}{{\mathcal L}}
\newcommand{\calL}{{\mathcal L}}
\newcommand{\cO}{{\mathcal O}}

\newcommand{\N}{\mathbb{N}}

\newcommand{\bsb}{{\boldsymbol b}}
\newcommand{\bsg}{{\boldsymbol g}}

\newcommand{\bsy}{{\boldsymbol y}}
\newcommand{\bsz}{{\boldsymbol z}}

\newcommand{\bnu}{{\boldsymbol \nu}}

\title[Multi-level Compressed Sensing discretization of high-dimensional parametric PDEs]{Multi-level Compressed Sensing Petrov-Galerkin discretization of high-dimensional parametric PDEs}
\thanks{
CS supported in part by European Research Council (ERC) through the grant AdG 247277, 
JLB and HR supported in part by the ERC grant StG 258926. 
JLB and HR would like to thank the Hausdorff Research Institute for Mathematics, 
University of Bonn, where parts of this work have been performed during 
the trimester program Mathematics of Signal Processing.
JLB thanks the support of the Clay Mathematics Institute for his visit to the CRM 
to attend the IRP Constructive Approximation and Harmonic Analysis 
where part of this work has been done.
}
\author{Jean-Luc Bouchot}
\address{Chair for Mathematics C (Analysis), RWTH Aachen University, Germany} \email{bouchot@mathc.rwth-aachen.de} 
\author{Holger Rauhut}
\address{Chair for Mathematics C (Analysis), RWTH Aachen University, Germany} \email{rauhut@mathc.rwth-aachen.de} 
\author{Christoph Schwab}
\address{Seminar for Applied Mathematics, ETH Z\"urich, Switzerland} \email{christoph.schwab@sam.math.ethz.ch}

\date{\today} 

\begin{document}

\maketitle


\begin{abstract}
We analyze a novel multi-level version
of a recently introduced compressed sensing (CS) 
Petrov-Galerkin (PG) method 
from 
[H. Rauhut and Ch. Schwab:
Compressive sensing {P}etrov-{G}alerkin approximation of
high-dimensional parametric operator equations,
Math. Comp. {\bf 304}(2017) 661--700]
for the solution of many-parametric partial differential equations. 
We propose to use multi-level PG discretizations, based on a hierarchy of nested 
finite dimensional subspaces, and to reconstruct parametric solutions 
at each level from level-dependent
random samples of the high-dimensional parameter space via CS methods such 
as weighted $\ell_1$-minimization.
For affine parametric, linear operator equations, we prove that
our approach allows to approximate the parametric solution
with (almost) optimal convergence order as specified by certain summability 
properties of the coefficient sequence in a general polynomial chaos expansion 
of the parametric solution and by the convergence order of the 
PG discretization in the physical variables.
The computations of the parameter samples of the PDE solution
is ``embarrasingly parallel'', as in Monte-Carlo Methods.
Contrary to other recent approaches, and as already noted in 
[A. Doostan and H. Owhadi:
A non-adapted sparse approximation of PDEs with stochastic inputs. 
JCP {\bf 230}(2011) 3015-3034]
the optimality of the computed approximations does not require a-priori assumptions
on ordering and structure of the index sets of 
the largest gpc coefficients 
(such as the ``downward closed'' property).
We prove that under certain assumptions work versus accuracy of the new algorithms 
is asymptotically equal to that of one PG solve 
for the corresponding nominal problem on the finest discretization level up to a constant. 
\end{abstract}

\section{Introduction}
\label{sec:Intro}
%
Motivated in particular by uncertainty quantification, 
the numerical solution of parametric operator equations 
has gained significant attention in recent years. 
In many cases, the underlying parameter 
space is high dimensional or even infinite dimensional
so that standard approximation methods
are subject to the curse of dimensionality, 
see e.g.\ \cite{Cohen15HighDPDEs,Chkifa14Breakingcurse}. 
Monte Carlo (\emph{MC}) sampling, however,
may be used in the context that the parametric model 
arises from a stochastic model and leads to a mean-square 
rate of $m^{-1/2}$ in terms of the number $m$ of sample evaluations, 
with constants that are independent of the 
parameter dimension.
The (dimension-independent) rate $1/2$ is not improvable in MC methods,
in general, and
the challenge consists in developing methods that achieve a faster 
convergence rate and at the same time alleviate or even overcome the curse of dimensionality.

A number of computational approaches have emerged in recent years 
towards this end. 
Among these are adaptive stochastic Galerkin methods, as developed
in \cite{EGSZ14_1230,EGSZ14_1045,G13_531},
reduced basis approaches (see, eg., \cite{BCDDPW11,BufMadPat12}),
adaptive Smolyak discretizations \cite{SS13_813,SS14_1117}, 
adaptive interpolation methods \cite{chkifa2014polyDownward} 
as well as sampling methods \cite{TangIacc2014}.
Adaptive Galerkin methods \cite{EGSZ14_1230,EGSZ14_1045,G13_531} are intrusive 
in the sense that they cannot simply reuse a solver developed for the corresponding problem 
with fixed parameter.
In contrast, the other above mentioned methods and algorithms are non-intrusive, 
but they rely on successive numerical solutions of the
operator equations for various parameter instances that are chosen 
based on suitable precomputations.
In contrast, (multilevel) Monte-Carlo (\emph{MLMC})\cite{MSS}, 
or Quasi-Monte Carlo approaches (\emph{QMC})\cite{DKLS16_1319}
compute expectations or statistical moments of the (random parametric) solution
via solutions for parameter instances chosen at random or ``quasi-random'', 
which allows to compute the ``parameter snapshot'' solutions in parallel. 

In this article, we build on a compressed sensing approach for 
numerically computing parametric solutions
developed and analyzed in \cite{Rauhut14CSPG,Bouchot15SLCSPG} 
(see also \cite{Doostan11sPDEsApprox,Peng14wl1PCE} for earlier work,
and \cite{ChkifDextTranWebs16} for recent developments)
and combine it with ideas originating from MLMC methods, 
see e.g.\ \cite{MLMCGilesActa,Barth11MLMCFE}. 
For Petrov-Galerkin (\emph{PG}) discretizations 
on a finite hierarchy of nested subspaces, 
ordered with respect to discretization levels,
the presently proposed method ``samples'', in a judicious fashion,
the parameter space and computes the corresponding PG approximations 
for random choices of the parameter vector.
As in MLMC-PG approaches, the number of such snapshot evaluations
decreases with increasing discretization level 
(corresponding to increasing refinement of the discretization).
In contrast to (ML)MC sampling, we employ
a CS technique based on weighted $\ell_1$-minimization 
\cite{Rauhut13wCS,Adck2015}
or iterative and greedy approaches (see for instance \cite{Jo2013wIHT,blda09,Foucart11HTP})
in order to reconstruct the coefficients of a generalized polynomial chaos
expansion of the difference of the parametric solution at two subsequent discretization levels. 
Finally, these differences are summed
together to obtain a PG approximation of the parametric solution at the finest level. 
One contribution of this paper is to show that the 
generalized polynomial chaos (GPC) expansion of the differences of 
PG approximations of the parametric solution is approximately sparse
by estimating the weighted $\ell_p$-norm for $0 < p < 1$
of the sequence of \Tscheb coefficients by a term that depends in 
a controlled way on the discretization level.
This fact makes the presently developed, multi-level version of the
compressive sensing approach feasible.   
We provide dimension-independent convergence rates 
which exceed $1/2$ under certain sparsity assumptions 
on the parametric solution family of the operator equation
and estimate the computational complexity for achieving such rates.
Similar to MLMC methods, 
the workload for approximating the parametric solution is
asymptotically the same as the one for computing one snapshot solution 
at the finest level up to a constant that depends only on smoothness parameters and $p \in (0,1)$.
%
However, in contrast to multilevel Monte Carlo, the convergence rates afforded by
our scheme are practically independent of the dimension and 
only limited by the solutions' sparsity; in particular, they may
significantly exceed $\mathcal{O}(m^{-1/2})$.

In mathematical terms, we consider linear, parametric operator equations of the 
generic form
\begin{equation}
\label{eq:operator}
A(\bsy)u(\bsy) = f.
\end{equation}
Here the parameter vector $\bsy \in U$ lies in a high-dimensional 
space $U$ making it challenging to computationally 
approximate the solution map $\bsy \mapsto u(\bsy)$, 
due to the mentioned
\emph{curse of dimensionality}, a notion going back to 
R.E. Bellman~\cite{BellmanCrseDim}, see~\cite{Cohen15HighDPDEs,Chkifa14Breakingcurse}
for its relevance in the present context.
Assuming that the parameter vector $\bsy = \(y_j\)_{j=1}^d$ 
takes values in finite intervals, we can consider, without loss of generality, 
$U = [-1,1]^d$, where the parameter set dimension $d$ may be finite or infinite. 

In our setting, the parametric family of operators $A(\bsy): \Xcal \to \Ycal'$ 
maps from a reflexive Banach space 
$\Xcal$ to the topological dual of, potentially, another reflexive Banach space $\Ycal$. 
A canonical example 
is the affine-parametric diffusion equation 
considered in \cite{cohen2010convergence,cohen2011analytic}
and in the single-level version of the present work~\cite{Rauhut14CSPG,Bouchot15SLCSPG}.
For a bounded Lipschitz domain $D \subset \R^n$ (one should think of $n=1,2,3$)
and a parametric diffusion coefficient $a(\cdot,\bsy)\in L^\infty(D)$
that depends affinely on a parameter vector $\bsy$, i.e.,
\be
\label{eq:aKL}
a(x,\bsy) = \bar{a}(x) + \sum_{j\geq 1} y_j \psi_j(x), \quad x\in D,
\ee
we consider the model parametric, second order 
divergence form elliptic Dirichlet problem
\be \label{eq:Diffusion}
A(\bsy) u := - \nabla \cdot(a(\cdot,\bsy) \nabla u) = f
\quad {\mbox{ in }} D,
\quad u|_{\partial D} = 0 .
\ee
The weak formulation of \eqref{eq:Diffusion}
in the Sobolev space $\bcX = \bcY :=H_0^1(D)$ reads:
Given $f\in \bcY'$, for every $\bsy\in U:= [-1,1]^\N$ 
find $u(\bsy)\in \bcX$ such that
\be\label{diff:eq:weak}
\int_D a(x,\bsy) \nabla u(x) \cdot \nabla v(x) \dif{x}
=
\int_D f(x) v(x) \dif{x}, \quad \mbox{ for all } v \in \bcY .
\ee
Eq.~\eqref{eq:Diffusion} is a particular example of an
\emph{affine-parametric operator equation} of the form
\begin{equation}
\label{eq:lineardependency}
A(\bsy) := A_0 + \sum_{j \geq 1} y_j A_j,\quad 
\bsy = (y_j)_{j\geq 1} \in U := [-1;1]^\Nbb,
\end{equation}
with $A_j := -\nabla \cdot (\psi_j \nabla), A_0 = -\nabla \cdot (\bar{a} \nabla)$.
In \eqref{eq:lineardependency}, 
the operator $A_0\in \calL(\bcX,\bcY')$ 
is traditionally referred to as \emph{nominal operator} or \emph{mean field}
while the operators $A_j\in \calL(\bcX,\bcY')$, for $j\geq 1$, 
are referred to as \emph{fluctuations}. 
For the parametric problem to be well-posed uniformly with respect to 
the parameter $\bsy\in U$, we assume that 
$\sum_{j\geq 1} \| A_j \|_{\calL(\bcX,\bcY')} < \infty$ in what follows.
Further assumptions required for the convergence and applicability of our approach
will be specified ahead.
Parametric expansions such as \eqref{eq:lineardependency} 
can be obtained e.g.\ by a Karhunen-Lo\`eve expansion
of random input data for divergence-form partial differential equations,
as explained in \cite{SchTodor06,cohen2011analytic}.

In order to ensure well-posedness of the parametric diffusion 
problem \eqref{eq:Diffusion} as in \cite{cohen2011analytic}
we require the {\em uniform ellipticity assumption}:
there exist constants $0<r \leq R <\infty$ such that
\be\label{UEA}
r \leq a(x,\bsy) \leq R, \qquad \mbox{ for almost all } x \in D, \mbox{ for all } \bsy \in U 
\;.
\ee
The Lax-Milgram Lemma ensures that for every $\bsy \in U$,
the weak formulation \eqref{diff:eq:weak} admits 
a unique solution $u(\cdot, \bsy) \in \bcX$ 
which satisfies the uniform a priori estimate
\[
\sup_{\bsy\in U} 
\| u(\bsy) \|_\Xcal \leq r^{-1} \|f\|_{\Ycal'}.
\]
Here and throughout the remainder, the term {\it ``uniform''} 
refers to uniform with respect to the parameter sequence $\bsy \in U$.

For the sake of simplicity we detail here only 
the approximation of functionals $\Gcal \in \Xcal'$ 
of solutions to the parametric operator equation~\eqref{eq:operator}, 
i.e., we are interested in the numerical approximation of
\[
F(\bsy) := \Gcal(u(\bsy)), \qquad \bsy \in U = [-1,1]^d,
\]
pointwise with respect to $\bsy$. 
We expect that our approach can be generalized to the
recovery of the vector-valued solution map $\bsy \mapsto u(\bsy)$, 
but we postpone this generalization to later contributions.
We are aiming at numerical schemes that are:
\begin{itemize}
\item Reliable: the convergence and accuracy should be verified and customizable;
\item Parallelizable: parallel sampling as in Monte-Carlo methods should be allowed,
with a convergence rate in terms of the number of samples
which (up to possibly logarithmic terms) equal the
best possible rate ensured by the compressibility of $F(\bsy)$, i.e., by weighted $\ell_p$-estimates
of the \Tscheb coefficients of $F$;
\item 
Non-intrusive: the approximation should use existing numerical solvers of the 
problem with fixed parameters, without any re-implementation of PDE solvers.
\end{itemize}
It is important to notice the difference to usual MC methods where the results obtained from random sampling usually hold in expectation. 
In contrast, our approach provides approximations that hold pointwise with respect to $\bsy$.
We estimate the coefficients of a tensorized \Tscheb expansion; whence
 only matrix-vector multiplications are required
in order to compute the solution $F(\bsy) = \mathcal{G}(u(\bsy))$
for any given parameter vector $\bsy =(y_j)_{j =1}^d$ up to a prescribed accuracy.
The computation scheme analyzed here differs from the single-level one introduced 
in~\cite{Rauhut14CSPG} in the sense that computing the approximation is done in a 
more efficient and computationally tractable manner. 
To this end, an unknown function $u(\bsy)$ is approximated by a telescopic sequence of 
so-called ``details''  at successively finer spatial resolutions:
$
u(\bsy) \approx \sum_{l=1}^L (u^l(\bsy) - u^{l-1}(\bsy))
$
where $u^l$ corresponds to a PG approximation on a discretization level $l$.
This is analogous to MLMC methods,
but is achieved here by \emph{compressive sensing} of the parameters 
with a \emph{level-dependent} number of parameter samples $\bsy^{(i)}$
on each discretization level in the physical domain.

We outline key ideas of the compressive sensing approach. 
We assume at our disposal a countable orthonormal basis 
$\(\varphi_\nu\)_{\nu\in \Lambda}$ of 
$L^2(U,\eta)$ with $\eta$ denoting a probability measure on 
the parameter set $U$ to be specified, 
and denote by $L^2(U,\eta;\Xcal)$ the Bochner space
of strongly measurable maps from $U$ to the (separable Hilbert) 
space $\Xcal$ containing solution instances, 
which are square integrable w.r.t.\ $\eta$.
We represent any function $u(\bsy)$ with values in $\Xcal$ 
as $u(\bsy) = \sum_{\nu \in \Lambda} \alpha_\nu \varphi_\nu (\bsy)$,
where 
${\bf \alpha} = \(\alpha_\nu\)_{j \in \Lambda}$ denotes the 
unique sequence of coefficients $\alpha_\nu \in \Xcal$.
Hence, in order to compute an approximation of the parametric solution for any $\bsy$ it suffices to
calculate an approximation of the coefficients $\alpha_\nu$.
For a new input parameter $\bsy$, 
one evaluates the basis functions $\varphi_\nu$ at $\bsy$ 
and forms a linear combination to recover a direct estimation of the solution.
Later on, we analyze the use of tensorized \Tscheb polynomials as orthonormal basis. 
The approximation is computed by evaluating the function at 
a few parameter points $\bsy^{(i)}$, $1 \leq i \leq m$, 
and solving the linear system $\bsg = \Phi {\bf \alpha}$,
where $\bsg = \(g_i\)_{i=1}^m = \(u(\bsy^{(i)})\)_{i=1}^m$ 
and where $\Phi$ corresponds to the sensing matrix
$\Phi \in \R^{m \times N}$ 
with entries $\Phi_{i,\nu} = \varphi_\nu(\bsy^{(i)})$, where $N$ corresponds to the number of basis functions
taken for the approximation.
However at this stage the coefficients $\alpha_\nu$ 
and the components $g_i$ are elements in $\Xcal$, and therefore, 
we first deal with the simpler case where a functional $\mathcal{G}$ 
(also known as the Quantity of Interest (\emph{QoI} for short) 
in the uncertainty quantification literature) 
is applied to the solution, resulting in
\[
b_i = \mathcal{G}(g_i) = \mathcal{G}\(u(\bsy^{(i)})\) = \sum_{\nu \in \Lambda}z_\nu \varphi_\nu(\bsy^{(i)}), \quad 
i = 1,\hdots, m, \qquad z_\nu = \mathcal{G}(\alpha_\nu), \quad \nu \in \Lambda
\;.
\]
We are particularly interested in the situation that the number $m$ of evaluations 
is smaller than the cardinality of $\Lambda$, so that
the linear system $b = \Phi z$ is underdetermined. 
Approximate sparsity of the coefficient sequence $(\alpha_\nu)$, and of $(z_\nu)$,
allows to apply techniques from compressive sensing such as (weighted) $\ell_1$-minimization 
or iterative hard thresholding (pursuit)
in order to recover $z$ accurately. 
In fact, approximate sparsity follows from the fact that $(\|\alpha_\nu\|_{\Xcal})$ and $(z_\nu)$ 
are contained in weighted $\ell_p(\mathcal{F})$-spaces, 
as shown in
\cite{Bachmayr2015sparsepoly,Chkifa14Breakingcurse,cohen2010convergence,
      cohen2011analytic,Rauhut14CSPG} 
and, for a related coefficient sequence, in this paper.

We expect that an approximation of the full solution $u(\bsy)$, $\bsy \in U$, 
taking values in the function space $\Xcal$, 
can be computed by a variant of our compressive sensing scheme. 
One may use ideas from joint/block sparsity 
in order to recover the sequence $(\alpha_\nu)$ with $\alpha_\nu \in \Xcal$ 
via mixed $\ell_1/\ell_2$-minimization, see e.g.~\cite{elmi09,elra10,Fornasier08jointSparse} 
(at least in the case that $\Xcal$ and $\Ycal$ are Hilbert spaces). 
However, we postpone a detailed analysis to a later contribution 
and restrict ourselves here to the simpler case of recovering the real-valued function
$\bsy \mapsto \mathcal{G}(u(\bsy))$.

The multi-level approximation scheme uses discretization levels $l=1,\hdots,L$, 
where the meshwidth at discretization level $l$ is $2^{-l} h_0$, 
so that the finest discretization is $h_L = 2^{-L} h_0$. 
With $n$ being the dimension of the domain $D$, we assume that the number of 
degrees of freedom at level $l$ scales like $\mathcal{O}(2^{nl})$, 
and 
we further assume 
\emph{available linear complexity, \changes{multigrid} solvers for the 
approximate solution of the discretized linear system of equations }
(uniformly with respect to the parameter $\bsy$)
resulting in computational costs per PG solution $u^l(\bsy^{(i)})$
that scales linearly in the number of degrees of freedom: 
$\mathcal{O}(2^{nl})$.

The presently proposed multi-level extension of the CS PG approach
from \cite{Rauhut14CSPG} proceeds analogous to MLMC
(see, e.g., \cite{Heinrich98MC} or \cite{MLMCGilesActa} and the references therein):
for parameter choices
$\{\bsy_l^{(i)}\}_{i=1,\hdots,m_l}$ on the discretization level $l$,
compute PG solutions $u^l(\bsy^{(i)}_l)$, $u^{l-1}(\bsy^{(i)}_l)$ 
at two consecutive discretization levels $l$ and $l-1$ (setting $u^{0} \equiv 0$). 
From the differences $\dfrak u^l(\bsy^{(i)}_l) = u^l(\bsy^{(i)}_l) - u^{l-1}(\bsy^{(i)}_l)$, 
we compute an approximation $\widehat{\dfrak u^l}(\bsy)$ 
via the single level compressive sensing approach of \cite{Rauhut14CSPG} 
for each $l = 1,\hdots, L$.
Finally, we combine the approximations at all levels 
similarly as in MLMC methods, i.e.,
$u^L_{\text{MLCS}}(\bsy) = \sum_{l=1}^L \widehat{\dfrak u^l}(\bsy)$, 
to obtain an approximation of the full parametric solution.
The main result of this paper consists of an analysis of
this method and provides, in its proof, 
a strategy on how to choose the number $m_l$ 
of parameter points at each level $l$. 
Its precise statement, Theorem~\ref{thm:mainThm}, 
is postponed to later in the exposition. 
To illustrate the type of results obtained here,
we state now a version of Theorem~\ref{thm:mainThm} in the particular case
of a linear, divergence form diffusion operator with affine dependence 
on the parameters (see Eqs.~\eqref{eq:aKL} and~\eqref{eq:Diffusion}). 
Ahead, we say that the weight sequence $\mathbf{v}$ is constant, 
if it is of the form $v_j = \beta$ for $j=1,\hdots,d$ 
for some $\beta > 1$ and $v_j = \infty$ for $j > d$, which corresponds
to the case that the expansion \eqref{eq:aKL} is finite (with $d$ terms). 
We say that $\mathbf{v}$ has polynomial growth if $v_j = c j^\alpha$, $j \in \N$, 
for some $c > 1$, $\alpha > 0$.
We refer to Section~\ref{sec:comp:costs} for details on the weight sequences.
\changes{Note that the order of the polynomial growth of the weight sequence does not affect the overall complexity (for a given target accuracy) of the method. It may however influence the multiplicative constants.}


%
\begin{atheorem}
\label{thm:simpleMainThm}
Let $L \in \Nbb$ and $\gamma \in (0,1)$. 
Consider the diffusion equation~\eqref{eq:Diffusion} with 
affine parametric coefficient~\eqref{eq:aKL}, forcing term 
$f \in H^{-1+t}(D)$ and functional $\Gcal \in H^{-1+t'}(D)$, 
with the respective smoothness parameters $t,t' \geq 0$. 
Assume that 
\changes{~\eqref{eq:aKL} holds with }$\bar{a}\in W^{t,\infty}(D)$ and that
the fluctuations fulfill the weighted 
$p$-summability\footnote{To ease the presentation, here and throughout the paper, 
we have not highlighted the dependence of the summability parameter $p$ 
on the regularity $t$ of the right-hand-side $f$.
It should be noted that the compressibility of the gpc expansion, 
the choice of the weight sequence, the number of samples per 
level all depend on the regularity of the data $\bar{a}$, $\psi_j$, $D$ and $f$.}
\begin{equation*}
\sum_{j \geq 1}\|\psi_j\|_{W^{t,\infty}(D)}^{p} 
v_{j}^{2-p} < \infty, 
\end{equation*}
for a sequence 
$\vbf = (v_j)_{j \geq 1}$
of weights 
as well as the following 
stronger, \emph{weighted} version of the 
Uniform Ellipticity Assumption~\eqref{UEA}: 
there exists $0 < r \leq R < \infty$ such that
\begin{equation}
\label{eq:UEAstrong}
 \sum_{j \geq 1}v_j^{(2-p)/p}|\psi_j(x)| \leq \min\{ \bar{a}(x)-r, R-\bar{a}(x) \}, 
\quad \text{for all } x \in D.
\end{equation}
%
With probability at least $1-\gamma$, 
the function $F(\bsy) := \Gcal(u(\bsy))$, $\bsy \in U$, 
can be approximated by $L$ (weigthed) sparse approximations (typically via weighted $\ell^1$ minimization)
based on a sequence of Galerkin projections into spaces of piecewise polynomials
    on regular, simplicial triangulations of meshwidth $h_\ell = 2^{-\ell} h_0$
from 
\begin{align*}
m_l &\gtrsim \max\{s_l \log^3(s_l) \log(N_l), \log(L/\gamma) \}
\end{align*}
%
solution evaluations at discretization level $l$ for $l = 1,\hdots,L$, 
where $s_l \asymp 2^{(L - l)(t+t')p/(1-p)}$, and $N_l$ is the size of the (level-dependent) 
active set $\Gamma_l$ of tensorized \Tscheb polynomials.
The resulting approximation $F^\#$ satisfies 
\begin{align*}
 \|F - F^\#\|_\infty &\leq C_p \|f\|_{H^{-1+t}(D)}\|\Gcal\|_{H^{-1+t'}(D)}L 2^{-(t+t')L} h_0^{t+t'}\\ 
 \|F - F^\#\|_2 &\leq C_p'\|f\|_{H^{-1+t}(D)}\|\Gcal\|_{H^{-1+t'}(D)} 2^{-(t+t')L} h_0^{t+t'} 
\end{align*}
Under the assumption that the computational cost of a single solve at level $l$ 
scales linearly with respect to the number of degrees of freedom, i.e.,
is $\Ocal(2^{nl})$ (for an $n$-dimensional domain $D$),
this result is achieved with a total work for the computation of 
snapshot solutions that scales 
as 
$ \Ocal\left(\max\left\{ 2^{nL}, L^\beta 2^{L(t+t')p/(1-p)} \right\}\right)$, 
where $\beta = 4$ for constant weights $\vbf$ and 
$\beta = 5$ for polynomially growing weights. 
The constant hidden in the $\mathcal{O}$-notation includes a factor of $\log(d)$ 
in the case of constant weights.
\end{atheorem}
We note in passing that in what follows, the 
estimates of the overall computational work do not account 
for the numerical solution of the (weighted) sparse approximation required for 
the compressed sensing estimation of the mapping $F$. 
We justify this convention by the observation 
that the computational cost of  $\ell_1$-minimization
is often of lower order compared to the total cost of evaluating the PDE samples.
\changes{Section~\ref{ssec:compTime} validates empirically this claim.}

Our theorem shows that in the case of sufficiently strong summability, 
i.e., $\frac{(t+t')p}{1-p} < n$, at a total cost that scales as a constant
times a single PDE solve at the finest discretization level $L$,
the multilevel CSPG (\emph{MLCSPG}) strategy can approximate 
a fixed function $F$ \emph{for any} parameter vector $\bsy \in U$.
This is analogous to what is afforded by MLMC methods, 
but the present MLCSPG strategy allows to achieve
any convergence rate afforded by the gpc summability, 
and allows to approximate the full parametric dependence, 
while MLMC only yields expectations (or moments).
Moreover, in our case the computational work scales 
favorably with decreasing $p$, which corresponds
to better sparse approximation rates implied by the weighted $p$-summability 
of (norms of) polynomial chaos coefficients of the parametric solution. 
To be more precise, in the case of higher smoothness $t+t' > 0$, 
we obtain an approximation error that scales with $h_L^{t+t'}$.
With a small enough value of $p$, 
we may exploit smoothness in the physical domain (allowing $t+t'$ such that  $\frac{(t+t')p}{1-p} < n$) and 
balance approximation error for the PDE solves.
    In contrast, the computational work required by MLMC to achieve an 
    expected approximation error scaling as $h_L^t$
    grows proportionally to $2^{2tL}$ when $2t \geq n$ 
    (where $t$ corresponds to the smoothness of the solution in the physical domain), 
    see \cite[Theorem 5.7]{Barth11MLMCFE}, 
    and there is no parameter $p$ in MLMC whose tuning allows to avoid such growth.

Nevertheless, we note that $t$ and $p$ may not be tuned 
independently:
in many instances increased smoothness $t$ leads to a larger value 
of the summability parameter $p$.

We emphasize that the tools and results developed here 
do not require a particular structure on the support set of the best approximation. 
It is often the case (see e.g.~\cite{chkifa2014polyDownward,Migliorati14discreteL2}), 
that proofs and/or methods require 
the sets of active indices in $N$-term gpc approximations be \emph{downward closed},
their approximation properties then being, in particular, independent of 
    the polynomial system adopted for implementation.
In constrast, the presently proposed, compressed sensing based approach
can recover (with high probability) \emph{any support set of active
multi-degrees of tensorized \Tscheb polynomial approximations}
(only assuming very rough knowledge
 of its location as provided by weighted $\ell_p$-estimates 
 of polynomial chaos coefficients), 
yet still providing quasi-optimal rates of convergence.
Moreover, apart from the $\ell_1$-minimization part of the algorithm, 
all function evaluations can be done in parallel.
\changes{We would like to point out that while we do not impose a particular structure on the coefficients, this structure is embedded in the choice of the coefficient sequence $\omega$. 
In particular, with the tensor product structure used in Eq.~\eqref{eq:omegas}, one notices that the weights will favor indices $\nu$ which have active components $\nu_j > 0$ associated to smaller $v_j$.
Setting smaller $v_j$'s will lead to a larger search space for the active coefficients $\alpha_\nu$. 
If the $v_j$'s scale inversely with the norm of the operators $A_0^{-1}A_j$ and these operators are ordered in decreasing order of norm, then the choice \eqref{eq:omegas} of $\omega$ yields a downward close structure of the level-dependent sets $\Gamma_l$ of active multi-indices.
This downward closeness comes however as a consequence of the choice of weights, and not as a requirement for the method to work.
}

Theorem~\ref{thm:simpleMainThm} is a particular case 
of our main Theorem~\ref{thm:mainThm} 
which we prove in Section~\ref{sec:MLCSPG} after 
recalling some basics about Petrov-Galerkin approximations 
in Section~\ref{sec:ParOpEq} combined with compressed sensing 
techniques in Section~\ref{sec:Tscheb}.
Section~\ref{sec:implem} deals with pratical aspects such as
truncating the dimension of the parameter space. 
The paper is finally concluded by 
numerical experiments to illustrate the theory
in Section~\ref{sec:numResults}.
%
\section{Petrov-Galerkin approximations}
\label{sec:ParOpEq}
%
We deal with the pointwise numerical approximations of 
the countably-parametric operator equation Eq.~\eqref{eq:operator}.
Numerically accessing the parametric solution map $U\ni \bsy\mapsto u(\bsy)$
at a fixed parameter instance $\bsy\in U$ requires discretization
of Eq.~\eqref{eq:operator} also in ``physical space''.
To this end, we introduce two dense, one-parameter 
families of discretization spaces 
$\{\Xcal^h\}_{h > 0} \subset \Xcal$ and $\{\Ycal^h\}_{h > 0} \subset \Ycal$
of equal finite dimensions 
$N^h := {\rm dim}(\Xcal^h) = {\rm dim}(\Ycal^h)$
and assume that the parametric operator $A(\bsy)$ fulfills the 
discrete and uniform $\inf-\sup$ conditions:
there exists a $\mu > 0$ such that for any $h >0$ and $\bsy \in U$
\begin{equation} 
\label{eq:infsupA}
\left\{
\begin{array}{rl}
\inf_{0 \neq v^h \in \bcX^h} 
\sup_{0 \neq w^h \in \bcY^h} 
&\frac{\langle A(\bsy) v^h, w^h \rangle}{\|v^h\|_\bcX\|w^h\|_\bcY} \geq \mu > 0 
\\ 
\inf_{0 \neq w^h \in \bcY^h} 
\sup_{0 \neq v^h \in \bcX^h} 
&\frac{\langle A(\bsy) v^h, w^h \rangle}{\|v^h\|_\bcX\|w^h\|_\bcY} \geq \mu > 0\;.
\end{array}
\right.
\end{equation}
The PG projections are defined as the solution to the following 
weak variational problems:
\begin{equation}
\label{eq:galerkinProjections}
\text{Find } u^h(\bsy) := G^h(\bsy)(u(\bsy)), \text{ such that } 
\langle A(\bsy)u^h(\bsy), v^h\rangle = \langle f, v^h\rangle \quad \mbox{ for all } v^h \in \Ycal^h
\;.
\end{equation}
%

We recall the following classical result (see for example~\cite[Chapter 6]{Brezzi13book}).
\begin{aprop}
\label{prop:linGh}
Let $\Xcal^h$ and $\Ycal^h$ be discretization spaces for the PG method, 
such that the uniform discrete $\inf-\sup$ conditions~\eqref{eq:infsupA}
are fulfilled and assume that the bilinear operator 
$\Xcal \times \Ycal \ni (u,w) \mapsto \langle A(\bsy)u,w\rangle$ is continuous,
uniformly with respect to $\bsy\in U$.

Then the PG projections $G^h(\bsy): \Xcal \to \Xcal^h$ 
are well-defined linear operators, 
whose norms are uniformly bounded with respect to the 
parameters $\bsy$ and $h$, i.e.,
\begin{align}
\label{eq:priorBound}
\sup_{\bsy\in U} \sup_{h>0}
\|u^h(\bsy)\|_{\Xcal} &\leq \frac{1}{\mu}\|f\|_{\Ycal'}, 
\\
\label{eq:normProjection}
\sup_{\bsy \in U} \sup_{h > 0} \|G^h(\bsy)\|_{\Lcal(\Xcal)} &\leq \frac{C}{\mu}
\end{align}
The Galerkin projections are uniformly quasi-optimal:
for every $\bsy\in U$ we have the a-priori error bound
\begin{equation}
\label{eq:quasiOptimal}
\|u(\bsy) - u^h(\bsy)\|_\bcX
\leq 
\(1+\frac{C}{\mu}\)\operatorname{inf}_{v^h \in \bcX^h} \|u(\bsy) - v^h\|_\bcX
\;.
\end{equation}
\end{aprop}
As is classical in the theory of polynomial approximation 
 (see, e.g.~\cite{Davis63InterpApprox,Rivlin69FuncApprox}),
we use a holomorphic extension of the parametric operator family 
$A(\bsy)$ to complex parameter sequences $\bsz \in \Ocal \supset U$, where $\Ocal$ 
is some suitable subset of the complex plane. 
Here,
when dealing with extensions of operators and solutions to 
parameters taking values in the complex domain, 
we identify the function spaces $\Xcal$ and $\Ycal$ 
with their complexifications  $\Xcal\otimes \{ 1, \ii \}$ and  $\Ycal\otimes \{ 1, \ii \}$ 
for the sake of simplicity.
We require the parametric operator
$\Ocal \ni \bsz \mapsto A(\bsz)$ 
to be holormorphic with respect to any finite set of variables and to be boundedly invertible. 
Hereby, a Banach-space valued mapping $z \mapsto R(z) \in E$ of a single complex variable
is said to be holomorphic (in some open domain $\Ocal$) if 
\[
\lim_{h \to 0} \frac{R(z_0 + h) - R(z_0)}{h}
\]
exists in $E$ for any $z_0 \in \Ocal$, with $h\to 0$ understood in $\Cbb$.
Note that our assumption on $A$ requires it to be 
holomorphic with respect to any component $z_j$ of $\bsz$ 
independently. Joint holomorphy with respect to an arbitrary,
finite subset of variables $\bsz' = (z_j)_{j \in \Lambda}$ with $|\Lambda| < \infty$ 
of $\bsz$ then follows from Hartogs' theorem.
In the sequel, we will often assume that the open set $\Ocal$, 
on which $\bsz \mapsto A(\bsz)$ is holomorphic, contains the product of 
Bernstein ellipses $\Ecal_\rho = \bigotimes_{j \geq 1}\Ecal_{\rho_j}$ 
with $\Ecal_\sigma := \{ (z+z^{-1})/2, z \in \Cbb: |z| = \sigma \}$.
In the case of complex-parametric operators,
the bounded invertibility of 
$A(\bsz)$ is equivalent to the 
\emph{complex discrete $\inf-\sup$ conditions}: 
there exists a constant $\mu_\Cbb > 0$ 
such that for any $h >0$ and $\bsz \in \Ocal$ 
\begin{equation} 
\label{eq:infsupAcomplex}
\left\{
\begin{array}{rl}
\inf_{0 \neq v^h \in \bcX^h} 
\sup_{0 \neq w^h \in \bcY^h} 
&\Re \frac{\langle A(\bsz) v^h, w^h \rangle}{\|v^h\|_\bcX\|w^h\|_\bcY} \geq \mu_\Cbb > 0, 
\\ 
\inf_{0 \neq w^h \in \bcY^h} 
\sup_{0 \neq v^h \in \bcX^h} 
&\Re \frac{ \langle A(\bsz) v^h, w^h \rangle}{\|v^h\|_\bcX\|w^h\|_\bcY} \geq \mu_\Cbb > 0 .
\end{array}
\right.
\end{equation}

Approximation results on discretization spaces are usually combined with 
prior knowledge of the regularity of the data. 
For this, we assume that there exists a $0 < t \leq \overline{t}$ 
such that the parametric family $A(\bsy)\in \calL(\bcX,\bcY')$ is regular in given
smoothness scales $\{ \bcX_t \}_{t \geq 0}$, resp. $\{ \bcY_t \}_{t\geq 0}$, satisfying:
%
\be\label{eq:SmoothPrimal}
\bcX = \bcX_0 \supset \bcX_1 \supset ... \supset \bcX_t \;, 
\quad 
\bcY = \bcY_0 \supset \bcY_1 \supset ... \supset \bcY_t \;. 
\ee
Here, the smoothness index $t\geq 0$ denotes, for example,
a differentiation order in a scale of Sobolev or Besov spaces.
These spaces are defined by interpolation for non integer indices. 
We shall also require a corresponding scale on the dual side, with $\Xcal'_t := (\Xcal')_t$, and $\Ycal'_t := (\Ycal')_t$: 
\be\label{eq:SmoothDual}
\bcX' = \bcX'_0 \supset \bcX'_1 \supset ... \supset \bcX'_t \;, 
\quad 
\bcY' = \bcY'_0 \supset \bcY'_1 \supset ... \supset \bcY'_t 
\;. 
\ee
Note carefully that with this notation, 
$(\bcX_t)'$ {\bf generally differs} from $\bcX'_t$.
For example, in the case of the diffusion equation~\eqref{eq:Diffusion}, 
one may choose $\Xcal = H^1_0(D)$ and $\Xcal_{t} = H^{1+t}_0(D)$. 
In this case, $\Xcal' = H^{-1}(D)=(H^1_0(D))'$
and 
$\(\Xcal_{t}\)' = H^{-1-t}(D) \neq H^{-1+t}(D) = \(\Xcal'\)_t =:\Xcal_t'$.

A first statement of 
solution regularity in the scales \eqref{eq:SmoothPrimal}, \eqref{eq:SmoothDual}
takes the form of 
uniform bounded invertibility of the family of 
parametric operators $A(\bsy)$:
\begin{equation}
\label{eq:invXtYt}
A(\bsy) \in \Lcal\(\Xcal_t, \Ycal_t'\), \quad  \mbox{ for all } \bsy \in U, 
\quad \text{ and } 
\sup_{\bsy \in U} \|A(\bsy)^{-1}\|_{\Lcal(\Ycal_t',\Xcal_t)} < \infty.
\end{equation}

For the PG discretization, we assume at hand two one-parameter families 
$\{\bcX^h\}_{h>0}$ and $\{ \bcY^h\}_{h>0}$ of $\bcX$ and of $\bcY$, 
respectively, with finite, equal dimension:
$N^h = {\rm dim} (\bcX^h) = {\rm dim} ( \bcY^h ) <\infty$.
We assume furthermore that 
$\{\bcX^h\}_{h>0}$ and $\{ \bcY^h\}_{h>0}$ 
are dense in $\bcX$ and in $\bcY$, respectively.
Here the discretization parameter $h>0$ usually stands for 
the meshwidth in finite element discretizations of fixed polynomial
degree, on a quasiuniform triangulation of the physical bounded, 
polyhedral domain $D$.
We assume that these spaces admit
\emph{\changes{linear} approximation properties} in the smoothness scales\changes{\footnote{Note that it would be possible to include the case of higher order FEM by propagating the order $k$ of FEM in the remaining parts of the estimations. We have chosen not to derive these results here to ease the presentation.}}
\eqref{eq:SmoothPrimal}, \eqref{eq:SmoothDual}, 
\begin{equation}
\label{eq:approxProperty}
\begin{array}{c} 
\operatorname{inf}_{w^h \in \bcX^h} \|w - w^h\|_\bcX \leq C_t h^t \|w\|_{\bcX_t}, 
\quad \mbox{ for all } w \in \bcX_t,
\\
\operatorname{inf}_{v^h \in \bcY^h} \|v - v^h\|_{\bcY} \leq C_{t'} h^{t'} 
\|v\|_{\bcY_{t'}}, \quad \mbox{ for all } v \in \bcY_{t'}.
\end{array}
\end{equation}
\cs{
Such approximation properties hold, for example, 
for the model Dirichlet problem \eqref{eq:Diffusion} in polytopal domain
$D\subset \mathbb{R}^d$ and for Lagrangian Finite Elements of polynomial
degree $p\geq 1$ on quasiuniform,
regular simplicial triangulations of $D$ of meshwidth $h$ with the choices
$\bcX = \bcY = H^1_0(D)$ and $\bcX_t = \bcY_t = H^{1+t}(D)\cap H^1_0(D)$
for $0\leq t,t' < \min\{ t^\ast, p \}$. Here, $t^\ast > 0$ denotes a
limit on isotropic Sobolev regularity of the solution of \eqref{eq:Diffusion}
in $D$ which is due to several factors: a) smoothness of $\partial D$,
b) smoothness of $f$ and c) smoothness of the parametric coefficient
$\bsy\mapsto a(\cdot,\bsy)$.
Alternative choices (with possibly larger ranges of $t^\ast$) are 
weighted (Kondrat'ev) spaces $\bcX_t = \bcY_t$ and 
Lagrangean Finite Elements of polynomial
degree $p\geq 1$ on locally refined regular simplicial triangulations of $D$ of meshwidth $h$.
}

Together with the bounded invertibility of the family of operators $A(\bsy)$, it holds:
\begin{equation*}
\|u(\bsy) - u^h(\bsy)\|_\Xcal \myeq{\leq}{Eq.~\eqref{eq:quasiOptimal}} 
C\inf_{v^h \in \Xcal^h} \|u(\bsy)-v^h\|_{\Xcal} \myeq{\leq}{Eq.~\eqref{eq:approxProperty}}
c_t h^t\|u(\bsy)\|_{\Xcal_t} \myeq{\leq}{Eq.~\eqref{eq:invXtYt}} C_t h^t\|f\|_{\Ycal_t'}.
\end{equation*}
Here, the constant $C_t$ depends on a uniform 
bound on
the inverse of the parametric operator in the appropriate smoothness space: 
$\sup_{\bsy\in U} \|A(\bsy)^{-1}\|_{\Lcal(\Ycal_t',\Xcal_t)}$, 
and on the smoothness parameter $t$, but not on the discretization parameter $h$.

Moreover, 
as we confine the exposition to
functionals of solutions $F(\bsy) = \Gcal(u(\bsy))$ 
for some $\Gcal(\cdot)\in \Xcal'$,
we assume \emph{adjoint regularity}, i.e.,
there exists $t' \geq 0$, 
such that $\Gcal \in \Xcal_{t'}'$, and such that the 
parametric adjoint solution $w_\Gcal(\bsy) \in \Ycal'$ of the problem
\begin{equation}\label{eq:AdEq}
A(\bsy)^*w_\Gcal(\bsy) = \Gcal
\end{equation}
satisfies $w_\Gcal(\bsy)\in \Ycal'_{t'}$ uniformly with respect to $\bsy$:
\begin{equation}\label{eq:AdReg}
\sup_{\bsy\in U} \| w_\Gcal(\bsy) \|_{\Ycal'} \leq C \| \Gcal \|_{\Xcal'_{t'}} 
\;. 
\end{equation}
Under the adjoint regularity \eqref{eq:AdReg},
the uniform parametric discrete inf-sup condition \eqref{eq:infsupA}
and the approximation property \eqref{eq:approxProperty},
an Aubin-Nitsche duality argument as, e.g., in \cite{Kuo2012qmcRdmPDE},
implies superconvergence: 
for any $\bsy \in U$, with $F^h$ the functional applied to the 
parametric PG solution $u^h(\bsy)$ defined in \eqref{eq:galerkinProjections} 
on discretization spaces of parameter $h$,
\begin{equation}
\label{eq:errorSLFEM}
|F(\bsy) - F^h(\bsy)| 
\leq 
C_{t+t'}h^{t+t'}\|f\|_{\Ycal_t'}\|\Gcal\|_{\Xcal_{t'}'}
\;.
\end{equation}
%
\section{Single-Level Compressed Sensing Petrov-Galerkin approximations}
\label{sec:Tscheb}
The multi-level compressed sensing PG 
(MLCSPG) discretization is a generalization of the single-level algorithms and 
results developed in \cite{Rauhut14CSPG}. 
Analogous to MLMC path simulations (see e.g.~\cite{MLMCGilesActa} and the references there) 
or MLMC Finite Element discretizations (see e.g.~\cite{Barth11MLMCFE}) 
the MLCSPG method described here
considers a sampling scheme from \cite{Rauhut14CSPG} 
with a number of sampling points depending on the discretization level.

Such compressed sensing reconstruction techniques have already
shown promise in the context of numerical solutions of 
PDEs on high-dimensional parameter spaces: we refer, for example,
to~\cite{yan2012stochastic,Doostan11sPDEsApprox,Peng14wl1PCE,Rauhut14CSPG,Bouchot15SLCSPG}.
\changes{Note that these approaches differ from other compressed-sensing based approaches that are used for efficiently computing a single snapshot, see for instance~\cite{Brugiapaglia15corsing,Brugiapaglia17corsingADR}. 
In their work, the authors do not use weighted versions of compressed sensing and only use compressibility in the spatial domain.
In comparison, our work considers a compression in the parameter space, and a recovery using weighted compressed sensing.}

The key idea in the works~\cite{yan2012stochastic,Doostan11sPDEsApprox,Peng14wl1PCE,Rauhut14CSPG,Bouchot15SLCSPG} is to decompose the solution of Eq.~\eqref{eq:operator} 
via its (tensorized \Tscheb or Legendre) 
polynomial chaos expansion with respect to the parameter vector $\bsy$.
A strongly measurable mapping $u: U \to \bcX : \bsy \mapsto u(\bsy)$
which is square (Bochner-) integrable with respect to the \Tscheb measure
$\dif\eta$ over $U$ can be represented as a gpc expansion, i.e.,
\begin{equation}
\label{eq:chebpolychaos}
u(\bsy) = \sum_{\nu \in \Fcal} u_\nu T_\nu(\bsy), 
\end{equation}
where in this case the coefficients in this expansion 
are functions $u_\nu \in \Xcal$. 
Here $\Fcal := \{ \nu \in \Nbb_0^\Nbb: |\operatorname{supp}(\nu)| < \infty \}$ 
is the set of multi-indices with finite support. 
The tensorized Chebyshev polynomials are defined as 
\begin{equation}
\label{eq:TensTscheb}
T_\nu(\bsy) 
= \prod_{j=1}^\infty T_{\nu_j}(y_j) 
= \prod_{j \in \operatorname{supp}(\nu)}T_{\nu_j}(y_j), 
\quad \bsy \in U, \quad \nu \in \Fcal,
\end{equation}
with the univariate \Tscheb polynomials defined by
\begin{equation}
\label{eq:Tscheb}
T_j(t) = \sqrt{2}\cos\(j \arccos(t)\), \qquad \text{and } T_0(t) \equiv 1
\;.
\end{equation}
Defining the probability measure $\sigma$ on $[-1;1]$ as $\dif{\sigma}(t) := \frac{\dif{t}}{\pi \sqrt{1-t^2}}$, 
%
%
the univariate \Tscheb polynomials $T_j$ defined
in~\eqref{eq:Tscheb} form an orthonormal system in $L^2([-1,1];\sigma)$ 
in the sense that 
\begin{equation*}
\int_{-1}^1 T_k(t)T_l(t) \dif{\sigma}(t) = \delta_{k,l}, \quad k,l \in \Nbb_0
\;.
\end{equation*}
Similarly, with the product measure 
\begin{equation*}
\dif{\eta}(\bsy) 
:= \bigotimes_{j \geq 1} \dif{\sigma}(y_j) 
 = \bigotimes_{j \geq 1} \frac{\dif{y_j}}{\pi\sqrt{1-y_j^2}},
\end{equation*}
the tensorized Chebyshev polynomials \eqref{eq:TensTscheb} 
are orthonormal with respect to $\eta$ in the sense that
\begin{equation*}
\int \limits_{\bsy \in U} T_\mu(\bsy)T_\nu(\bsy) \dif{\eta}(\bsy) 
= 
\delta_{\mu,\nu}, \quad \mu,\nu \in \Fcal
\;.
\end{equation*}
A result proven in~\cite{Hansen13chebyshev} ensures the 
$\ell_p$ summability, for some $0 < p \leq 1$, of the polynomial 
chaos expansion~\eqref{eq:chebpolychaos} for the diffusion case, 
Eq.~\eqref{eq:Diffusion}: 
\begin{equation*}
 \left\|\(\|u_\nu\|_\Xcal\)_{\nu \in \Fcal}\right\|_p^p 
= 
\sum_{\nu \in \Fcal}\|u_\nu\|_\Xcal^p < \infty
\end{equation*}
under the uniform ellipticity assumption \eqref{UEA} and the condition that the sequence of infinity norms of the 
$\psi_j$ is itself $\ell_p$ summable:
\begin{equation*}
 \left\|\(\|\psi_j\|_\infty\)_{j \geq 1}\right\|_p^p = \sum_{j \geq 1}\|\psi_j\|_\infty^p < \infty.
\end{equation*}
Recent results by~\cite{Bachmayr2015sparsepoly} show that these conditions 
can be improved by considering pointwise convergence of the 
series $\sum_{j \geq 1}|\psi_j|$ instead of infinity norms 
in the whole domain $D$.
This takes advantage of the local structure of the basis elements 
$\psi_j$, e.g. when only few of them are overlapping, as is the case for wavelets. 
In particular, $\ell_p$ summability of Legendre and \Tscheb coefficients 
can be obtained when $\(\|\psi_j\|_\infty\)_{j} \in \ell_q$ for $q := 2p/(2-p)$
provided that the interiors of the supports of the $\psi_j$ do not overlap. 
%
The summability results from~\cite{Hansen13chebyshev} concerning \Tscheb expansions 
were extended to weighted $\ell_p$ estimates 
for the general parametric operator problem \eqref{eq:operator} with affine dependence as 
in~\eqref{eq:lineardependency}, in~\cite{Rauhut14CSPG}
%
%
under slightly stronger assumptions.
This result is particularly important for us as it ensures the recovery 
of the coefficients $u_\nu$ (or any functional thereof) via CS methods. 

The results on the approximation via an MLCSPG framework 
rely on the single-level results developed in~\cite{Rauhut13wCS,Rauhut14CSPG}, 
where functions are approximated via a weighted-sparse expansion in an appropriate basis. 
We review here the main ideas. 
Given a (finite) orthonormal system 
$\(\phi_\nu\)_{\nu \in \Lambda}$, with $|\Lambda| = N < \infty$ for $L^2(U,\eta)$ 
where $\eta$ is a probability measure, for any fixed function $f: U \to \Rbb$, 
there exists a unique sequence of coefficients 
$\fbf = \(f_\nu\)_{\nu \in \Lambda}$ such that 
\begin{equation}
\label{eq:expansionf}
 f(\bsy) = \sum_{\nu \in \Lambda} f_\nu \phi_\nu(\bsy), ~\forall y \in U\;.
\end{equation}
We define an $\ell_p$ norm associated with this expansion as 
$\vertiii{f}_p := \|\fbf\|_p$.

In particular, a function $f$ is said to be sparse (or compressible) 
if its sequence of coefficients in expansion~\eqref{eq:expansionf} 
is sparse (or compressible) itself. 
Our goal is to recover this said sequence from seemingly few evaluations 
of the function $f$ at certain (here random) sample points $\bsy^{(i)}$, 
for $1 \leq i \leq m$. 
This can be done by CS methods:
after introducing the sensing matrix $\Phi$ as $\Phi_{i,j} := \phi_j(\bsy^{(i)})$ 
and letting $g_i := f(\bsy^{(i)})$, it holds
\begin{equation*}
 \bsg = \Phi \fbf\;.
\end{equation*}
Hence, assuming that the expansion is sparse, and the number of samples 
rather small, we are dealing with the by-now classical problem of recovering 
a sparse vector from few linear measurements, 
by solving, for instance, the convex program
\begin{equation}
\label{eq:BP}
\min_{\bsz \in \Rbb^N} \|\bsz\|_1, \quad \text{subject to } \Phi\bsz = {\bsg}.
\end{equation}
%
In our context, it is beneficial to use a \emph{weighted framework} which has been developed recently in~\cite{Rauhut13wCS}.
Introducing a sequence of positive weights 
$\(\omega_{\nu}\)_{\nu \in \Lambda}$ with $|\omega_\nu| \geq 1$ for all $\nu$, 
a weighted $\ell_p$ (quasi-)norm (henceforth indexed $\ell_{\omega,p}$ when appropriate) 
can be defined as 
\begin{equation*}
\|\fbf\|_{\omega,p}^p := \sum_{\nu \in \Lambda} \omega_\nu^{2-p}|f_\nu|^{p}, \quad 0< p \leq 2.
\end{equation*}
In particular, it holds 
$\|\fbf\|_{\omega,2} = \|\fbf\|_2$ and 
$\|\fbf\|_{\omega,1} = \|\fbf \odot {\bf \omega}\|_{1}$, 
where $\odot$ defines the pointwise multiplication. 
Moreover, choosing the constant weight $\omega_\nu = 1$ yields the original 
definitions of $\ell_p$ norms.
%
Formally letting $p \downarrow 0$ motivates the introduction of
    the \emph{weighted sparsity} measure
\begin{equation*}
\|\fbf\|_{\omega,0} := \sum_{\nu \in \Lambda, f_\nu \neq 0} \omega_\nu^2
\;.
\end{equation*}
A vector $\xbf$ is therefore called \emph{weighted $s$-sparse} 
(with respect to a weight sequence $\omega$) if $\|\xbf\|_{\omega,0} \leq s$. 
We may therefore define the error of best weighted $s$-term approximation 
as
\begin{equation*}
\sigma_{\omega,s}(f) = \sigma_{\omega,s}(\fbf) 
:= \inf_{\bsz: \|\bsz\|_{\omega,0} \leq s} \|\fbf-\bsz\|_{\omega,p}.
\end{equation*}

With these weighted error measures at hand,
the Basis Pursuit problem~\eqref{eq:BP} can be 
generalized to include a-priori information encoded in the 
sequences $\omega$ of weights, as 

\begin{equation}
\label{eq:wBP}
\min_{\bsz \in \Rbb^N} \|\bsz\|_{\omega,1}, \quad \text{subject to } \Phi\bsz = {\bsg}.
\end{equation}
More details on such weighted spaces and weighted sparse approximations 
can be found in~\cite{Rauhut13wCS} where the following fundamental result 
is also proved.
\begin{atheorem}
\label{thm:wl1approx}
Suppose $\(\phi_\nu\)_{\nu \in \Lambda}$ is a finite orthonormal system with 
$|\Lambda| = N < \infty$ and that weights $\omega_\nu \geq \|\phi_\nu\|_{\infty}$ are given. 
For a (weighted) sparsity $s \geq 2\|{\bf \omega}\|_{\infty}^2$, draw 
\begin{equation}
\label{eq:mwl1approx}
\changes{m \geq Cs \max\{\log^3(s)\log(N), \log(1/\gamma)\}}
\end{equation}
sample points $\bsy^{(i)}$ at random, $1 \leq i \leq m$,
according to the orthonomalization measure $\eta$.
The constant $C>0$ in \eqref{eq:mwl1approx} is universal, i.e., independent of all other quantities including $s$, $m$ and $N$.

Then, with probability at least \changes{$1-\gamma$}, 
any function 
$f = \sum f_\nu \phi_\nu$ can be approximated by the function 
$\ucs{f} := \sum \ucs{f}_\nu\phi_\nu$, where $\ucs{\fbf}$ 
is the solution to the weighted basis pursuit problem~\eqref{eq:wBP}.
The approximation holds in the following sense:
$$
\|f-\ucs{f}\|_{\infty} \leq \vertiii{f-\ucs{f}}_{\omega,1} 
\leq c_1 \sigma_s(f)_{\omega,1}, 
\quad \text{ and } \quad 
\|f-\ucs{f}\|_2 \leq d_1 \sigma_s(f)_{\omega,1} / \sqrt{s}
\;.
$$
\end{atheorem}

In particular, using the weighted Stechkin inequality from~\cite{Rauhut13wCS}
\begin{equation}
\label{eq:stechkin}
\sigma_s(\fbf)_{\omega,q} 
\leq 
\(s-\|\omega\|_\infty^2\)^{1/q-1/p}\|\fbf\|_{\omega,p}, 
\quad p < q \leq 2, \quad \|\omega\|_\infty^2 < s
\;,
\end{equation}
we obtain \cs{that for given summability exponent $0 < p < 1$, 
there exists a constant $C>0$ independent of $s$ 
such that
}
\begin{align}
\label{eq:approxSampPts}
\|f-\ucs{f}\|_\infty &\leq \cs{C} s^{1-1/p} \|\fbf\|_{\omega,p}, \quad
 \|f-\ucs{f}\|_{2} 		%
\leq \cs{C} s^{1/2-1/p} \|\fbf\|_{\omega,p}.
\end{align}

Choosing $m \asymp s \ln^3(s) \log(N)$ relates the reconstruction error and the number $m$ of samples as
\[
\|f-\ucs{f}\|_\infty \leq c \left(\frac{\log^3(m) \log(N)}{m}\right)^{1/p-1} \|\fbf\|_{\omega,p}, \quad
 \|f-\ucs{f}\|_{2} 		%
\leq d  \left(\frac{\log^3(m) \log(N)}{m}\right)^{1/p-1/2} \|\fbf\|_{\omega,p}.
\]

\begin{rmk}
\label{rmk:nbSamples}
Recent works~\cite{Bourgain2014RIP,Haviv15RIPFT} on restricted isometry constants for subsampled 
Fourier matrices suggest that the factor $\log^3(s)$ in~\eqref{eq:mwl1approx} can be reduced to $\log^2(s)$.
\end{rmk}

%
\section{Multi-level Compressed Sensing Petrov-Galerkin approximations}
\label{sec:MLCSPG}
%
\subsection{A multi-level framework}
\label{sec:MLSetting}
We extend the foregoing CS methods to sweeping the parameter domain
to multi-level (``ML'' for short) discretizations of the parametric problems, 
in the spirit of the ML MC methods for numerical treatment of operator 
equations with random inputs as 
developed in~\cite{Heinrich98MC,MLMCGilesActa,Barth11MLMCFE}. 
There,
the solution of the parametric operator equation~\eqref{eq:operator} 
is approximated on a sequence of partitions of the physical domain $D$ of widths
$\{h_l\}_{l=1}^L$ for a prescribed, maximal refinement level $L \in \Nbb$. 
To simplify the exposition, we assume dyadic refinement, i.e. 
$h_{l+1} = h_l/2 = 2^{-l}h_0$ for a given, small enough, 
initial resolution $h_0 > 0$.

For a given parameter sequence $\bsy$, we may write the Galerkin
projection $u^L(\bsy) \in \Xcal^{h_L}$ of $u(\bsy)$ as
\begin{equation}
\label{eq:mldecomposition}
u^L(\bsy) = \sum_{l = 1}^L u^{l}(\bsy) - u^{l-1}(\bsy),
\end{equation}
where we define $u^0(\bsy) \equiv 0$ 
(note that we will equivalently parametrize the approximations and spaces by $l$ or $h_l$). 
The idea behind our MLCSPG approach
is to estimate every difference between consecutive levels of approximation (the \emph{details}) via a single level CSPG as presented above. 
For the remaining, we let 
\begin{equation}\label{def:du}
\deltau{u}^{l}(\bsy) := u^l(\bsy) - u^{l-1}(\bsy), \quad 1 \leq l \leq L,
\end{equation}
denote the difference between two scales of approximation. 

As already outlined in the introduction, our method produces pointwise
numerical approximations ${\ucs{\deltau{u^l}}}(\bsy)$
of $\deltau{u}^{l}(\bsy)$ via a (single level) CSPG method. 
For each level $l$, we choose a number $m_l$ of parameter vectors $\bsy_l^{(1)},\hdots,\bsy_l^{(m_l)}$, compute
the PG approximations $u^{l}(\bsy_l^{(i)})$ and $u^{l-1}(\bsy_{l}^{(i)})$ by solving the corresponding 
finite dimensional linear systems, and form the samples 
$\deltau{u}^{l}(\bsy_l^{(i)}) = u^l(\bsy_l^{(i)}) -  u^{l-1}(\bsy_l^{(i)})$, $i = 1,\hdots,m_l$. 
From these samples, one approximates the coefficients in the tensorized Chebyshev expansion of $\deltau{u}^{l}$
via weighted $\ell_1$-minimization (or any sparse recovery method).
This yields approximations $\ucs{\deltau{u^l}}(\bsy)$, $\ell=1,\hdots,L$, and 
\[
u^L_{\text{MLCS}}(\bsy) := \sum_{l=1}^L \ucs{\deltau{u^l}}(\bsy)
\] 
then provides an approximation of the targeted parametric solution $u = u(\bsy)$.
%
The convergence of our MLCSPG framework can be estimated
via the triangle inequality, 
\begin{equation}
 \label{eq:errorToAnalyzeU}
 \|u(\bsy) - {u^L}_{\text{MLCS}}(\bsy)\|_\Xcal 
  \leq \|u(\bsy) - u^L(\bsy)\|_\Xcal 
  + \sum_{l=1}^L \left\| \dfrak u^{l}(\bsy) - \ucs{\deltau{u}^l}(\bsy)\right\|_\Xcal.
\end{equation}
For simplicity, we constrain our considerations to a functional $\Gcal \in \Xcal'$ 
applied to the solution, leading to the real-valued 
QoI $F(\bsy) = \Gcal(u(\bsy))$ to be approximated. 
The above considerations apply verbatim when replacing $u(\bsy)$ by $F(\bsy)$, 
and $\deltau{u}^{l}$ by $\dfrak F^l$, the levelwise PG approximation, 
and $\ucs{\deltau{u^l}}(\bsy)$ by $\ucs{\deltau{F^l}}(\bsy)$.
The triangle inequality leads to the error estimate
\begin{equation}
 \label{eq:errorToAnalyze}
 \left| F(\bsy) - {F^L}_{\text{MLCS}}(\bsy)\right| 
\leq 
|F(\bsy) - F^L(\bsy)| + \sum_{l=1}^L \left| \dfrak F^{l}(\bsy) - \ucs{\deltau{F}^l}(\bsy)\right|.
\end{equation}
The first term on the right hand side of Eq.~\eqref{eq:errorToAnalyzeU}
can be estimated using the uniform parametric regularity \eqref{eq:invXtYt},
the uniform parametric inf-sup condition \eqref{eq:infsupA} 
and the approximation property \eqref{eq:approxProperty}:
for a regularity parameter $0 < t \leq \overline{t}$ 
of the data $f$, 
%
\begin{equation*}
\|u(\bsy) - u^L(\bsy)\|_\Xcal \leq C_t h_L^t \|f\|_{\Ycal_t'}.
\end{equation*}
%
Passing to the functional $\Gcal \in \Xcal_{t'}'$,
we obtain a superconvergence error bound for the Petrov-Galerkin-Finite Element Method (PG-FEM)
via a classical Aubin-Nitsche duality argument~\cite{Kuo2012qmcRdmPDE}:
\begin{equation*}
|F(\bsy) - F^L(\bsy)| \leq C_{t+t'}h_L^{t+t'}\|f\|_{\Ycal_t'}\|\Gcal\|_{\Xcal_{t'}'}.
\end{equation*}
Our goal is to verify that the single-level result,
Theorem~\ref{thm:wl1approx}, applies to all levels $l=1,\hdots L$, 
and to obtain error bounds similar to the one in Eq.~\eqref{eq:errorSLFEM}. 
We consider the \Tscheb expansions of the differences, 
\begin{align}
\label{eq:expUdl} 
\mathfrak{d}u^{l}(\bsy) &= \sum_{\nu \in \Fcal} \mathfrak{d} u_\nu^{l} T_\nu(\bsy),
\\
\label{eq:expFdl} 
\dfrak F^{l}(\bsy) &= \sum_{\nu \in \Fcal} \dfrak F_\nu^{l} T_\nu(\bsy).
\end{align}
Assuming summability of the expansion in $\ell_{\omega,p}(\Fcal)$, 
we can apply Theorem~\ref{thm:wl1approx} with a number of samples 
$m_l \gtrsim s_l \log(s_l)^3 \log(N_l)$, for suitable choices of $s_1,\hdots,s_L$,
and in particular we can use the error 
estimate \eqref{eq:approxSampPts} in terms of the (weighted) sparsity 
$s_l$ for each level of approximation.
This results in the bound
\begin{equation}
 \label{eq:approxSampPtsTscheb}
 |\dfrak F^{l}(\bsy) - \ucs{\deltau{F}^l}(\bsy)| 
  \leq 
  C \left\|\(\dfrak F^{l}_{\nu}\)_{\nu}\right\|_{\omega,p} s_l^{1-1/p}, 
 \quad \text{for all } 1 \leq l \leq L,
\end{equation}
where $C > 0$ is a universal constant (independent of $s_l$, $\bsy$, $l$).
Theorem~\ref{thm:wl1approx} applies only to finite orthonormal systems.
Thus, for each $l=1,2,...,L$,
the countably infinite index set $\Fcal$ has to be truncated to a finite, 
but possibly large, subset $\Gamma_l$ of $N_l := |\Gamma_l| < \infty$ 
many indices of the relevant (few) essential \Tscheb coefficients 
in the parametric solution's gpc expansion.
We describe a strategy for selecting the index sets
$\Gamma_l$ depending on $s_l$ in Section~\ref{ssec:truncation}.
A good choice for the $s_l$ turns out to be $s_l \asymp 2^{(L-l)(t+t')p/(1-p)}$, 
as will be derived ahead.

Finally, summing up the contributions from all discretization levels
and drawing 
$$m_l \gtrsim  2^{(L-l)(t+t')p/(1-p)}(L-l)^3 \log({N_l})$$ sample points per level 
will imply the error bounds in Theorem~\ref{thm:mainThm}. 
The choice of this number of sampling points is 
justified in Section~\ref{sec:convergence} 
and by the following result, whose proof is the purpose of the next section.

\begin{atheorem}
\label{thm:summability}
Let $\{A(\bsy): \bsy \in U\}$ be a parametric family of operators as defined
in \eqref{eq:lineardependency}.
Assume that the operator $A_0$ is $\inf-\sup$ stable. 
For $B_j := A_0^{-1}A_j$ and for $0\leq t \leq \bar{t}$, 
introduce the sequence 
\begin{equation}\label{eq:Bj}
\bsb_t 
:= \(b_{t,j}\)_{j\geq1} \;\;\text{with}\;\; 
b_{t,j} := \|B_j\|_{\Lcal(\Xcal_t)} = \| B_j^* \|_{\Lcal(\Ycal_t)}.
\end{equation}
Let $\vbf := \(v_j\)_{j \geq 1}$ be a sequence of weights with $v_j \geq 1$ such that, 
for some $p < 1$,
\begin{align} \label{eq:sumbjvjtilda} 
\sum_{j \geq 1} b_{t,j} v_j^{(2-p)/p} &\leq \kappa_{\vbf,p} < 1, \text{ and }
\\
\label{eq:psumbjvj} 
\sum_{j \geq 1} b_{t,j}^p v_j^{2-p} &< \infty.
\end{align}
Let $\rho$ be a \changes{$\bsb_0-\delta$-admissible} sequence of polyradii, 
with $\delta = (1-\kappa_{\vbf,p})/2$, i.e., such that 
\begin{equation} \label{eq:admissibility}
\sum_{j \geq 1}(\rho_j -1)b_{0,j} \leq \delta.
\end{equation}
Then the family of operators $A(\bsy)$ is uniformly $\inf-\sup$ stable.
Assume in addition that 
$A_j \in \calL(\Xcal_t,\Ycal_t')$, $j \geq 0$, 
are defined on the scale of smoothness spaces $\Xcal_t$ 
and that the approximation property~\eqref{eq:approxProperty} holds.
Assume moreover that 
$A_0: \Xcal_t \to \Ycal_t'$ is boundedly invertible 
and 
that the sequence $\bsb_t$ is small, 
and that the polyradius $\rho$ is \changes{$\bsb_t-\delta_t$-admissible}, i.e.
\begin{align}
\label{eq:sumbjt} \sum_{j \geq 1} b_{t,j} \leq \kappa_{t} &< 1, 
\\
\label{eq:deltatAdmissibility}\sum_{j \geq 1} (\rho_j - 1) b_{t,j} &\leq \delta_t,
\end{align}
for $\delta_t < 1-\kappa_t$\;.

Then the affine-parametric family of operators 
$\{A(\bsy):\; \bsy\in U\}$ 
is uniformly boundedly invertible in $\Lcal(\Xcal_t, \Ycal_t')$,
and there hold bounds on the \Tscheb gpc coefficients
\begin{equation}\label{eq:ChdltX}
\|\dfrak u_\nu^{l}\|_{\Xcal} \leq Ch_l^t\|f\|_{\Ycal_t'}\rho^{-\nu}\;,
\quad \mbox{and} \quad 
|\dfrak F_\nu^{l}| \leq Ch_l^{t+t'}\|\Gcal\|_{\Xcal_{t'}'}\|f\|_{\Ycal_t'}\rho^{-\nu}
\quad \mbox{for all}\quad \nu \in \Fcal.
\end{equation}
Moreover, 
for each $\nu \in \Fcal$, there exists a $\delta$-admissible sequence 
$\rho = \rho(\nu)$ satisfying  \eqref{eq:deltatAdmissibility} such that
the sequence with components 
$\rho(\nu)^{-\nu} = \prod_{j \geq 1} \rho(\nu)_j^{-\nu_j}$, $\nu \in \Fcal$, 
satisfies 
$\(\rho(\nu)^{-\nu}\)_{\nu \in \Fcal} \in \ell_{\omega,p}$, 
where 
\begin{equation}\label{def:weight:omega}
\omega_\nu := \theta^{\|\nu\|_{0}}\vbf^{\nu} = \theta^{\|\nu\|_0} \prod_{j\geq 1} v_j^{\nu_j}, 
\quad \nu \in \Fcal\;.
\end{equation}
\end{atheorem}
We want to stress once again that the result presented above is written without the explicit dependence of the weight sequence $\vbf$ on the regularity parameter $t$.
Moreover, we note that the conditions \eqref{eq:Bj} - \eqref{eq:deltatAdmissibility}
are, for $t>0$, strictly stronger than the summability conditions
which were required 
in the single-level PG analysis in \cite{Rauhut14CSPG}.

\subsection{Summability of the \Tscheb expansions}
\label{sec:SummTscheb}
This section provides the proof of the core result of the present paper,
Theorem \ref{thm:summability}.
We show that under general assumptions, 
the parametric solution's sequence of \Tscheb coefficients
$\(\dfrak F_\nu^{l}\)_{\nu \in \Fcal} \in \ell_{p,\omega}$, 
and in particular that the following a priori estimate holds:
\begin{equation}
\label{eq:priorBoundFnu}
\left\|\(\dfrak F_\nu^{l}\)_{\nu \in \Fcal}\right\|_{\omega,p} 
  \leq Ch_l^{t+t'}\|f\|_{\Ycal_t'}\|\Gcal\|_{\Xcal_{t'}'} \|(\rho(\nu)^{-\nu})_{\nu \in \Fcal}\|_{\omega,p}.
\end{equation}
The main novel point of this estimate is the scaling of the right hand side with 
$h_l^{t+t'}$.
The proof of this assertion is structured in three main steps, 
analogously to~\cite{cohen2011analytic,Rauhut14CSPG}.
First we show that the difference between levels is holomorphic in polydiscs. 
Then, this holomorphy is used to bound the norm of any \Tscheb coefficient. 
This norm depends on a sequence of radii of holomorphy $\rho$. 
Finally, we construct a sequence of radii and weights such that the 
sequence of coefficients is $\ell_{\omega,p}$ summable. 
\subsubsection{Holomorphy}
\label{sec:holomorphy}
This first part shows that, 
under some uniform invertibility assumption of the family of (complexified) operators $A(\bsz)$
(which are satisfied in particular for the affine-parametric family considered here),
the solutions are holomorphic with respect to any finite set of variables. 
This then allows to use Cauchy's integral formula to estimate the norm 
of the \Tscheb coefficients.
\begin{atheorem} \label{thm:CplxInfSup}
For some $\Ocal \subset \C^\N$ with $\Ocal \supset U$, 
assume that the complex $\inf-\sup$ conditions~\eqref{eq:infsupAcomplex} hold 
with constant $\mu_\Cbb$ uniformly for $\bsz \in \Ocal$. 
If 
the solution map $\Ocal \ni \bsz \to u(\bsz) \in \Xcal$ 
is holomorphic with respect to any finite set of parameters, then
\begin{enumerate} %
\item 
for any level $l$ of PG discretization  
(corresponding to the discretization parameter $h_l = 2^{-l}h_0$ 
 for a given $h_0 > 0$ sufficiently small),
the parametric Galerkin projections $\Ocal \ni \bsz \to u^l(\bsz)\in \Xcal^l$ 
are holomorphic with respect to any finite subset of the sequence $\bsz\in \Ocal$,
with domains of holomorphy whose size is independent of $l$, i.e. of
the discretization parameter $h_l$,
\item 
the Petrov-Galerkin projections are quasi-optimal, uniformly with respect
to the level of approximation $l$ and the vector of (complex) 
parameters $\bsz \in \Ocal$:
\begin{equation*}
\|u(\bsz) - u^l(\bsz)\|_\Xcal 
\leq 
\(1 + \frac{C}{\mu_\Cbb}\)\inf_{v^l \in \Xcal^l} \|u(\bsz)-v^l\|_{\Xcal}
\;.
\end{equation*}
\end{enumerate}
\end{atheorem}

\begin{proof}
The holomorphy follows from the linearity of the 
PG approximation as stated in Proposition~\ref{prop:linGh}. 
The quasi optimality is obtained in the same way as in the real case.
\end{proof}
%
The next corollary which uses the notation \eqref{def:du} follows directly. 
\begin{acorol}
\label{cor:normDiff}
Under the conditions above, if in addition the approximation 
property of the discretization spaces holds for complex parameters $\bsz \in \Ocal$, 
then for any two consecutive discretization 
levels $l$ and $l+1$, $l \geq 0$, the mappings 
$\Ocal \ni \bsz \mapsto \deltau{u}^l(\bsz) \in \Xcal$ 
are holomorphic with respect to any finite set of variables and 
satisfy the uniform bound 
\begin{equation*}
\sup_{\bsz \in \Ocal} \|\deltau{u}^l(\bsz)\|_{\Xcal} 
\leq 
C_{t,\mu_\Cbb}' h_l^t\sup_{\bsz \in \Ocal} \|u(\bsz)\|_{\Xcal_t}.
\end{equation*}
\end{acorol}

\begin{proof}
The statement is a consequence of the previous results and the triangle inequality:
\begin{align*}
\sup_{\bsz \in \Ocal} \|\deltau{u}^l(\bsz)\|_{\Xcal} 
&~ \leq \(1 + \frac{C}{\mu_\Cbb}\)\sup_{\bsz \in \Ocal} 
\(
\inf_{v^{l} \in \Xcal^{l}}\|u(\bsz)-v^{l}\|_{\Xcal} 
      + \inf_{v^{l-1} \in \Xcal^{l-1}}\|u(\bsz)-v^{l-1}\|_{\Xcal}\)
\\ 
& \leq \(1 + \frac{C}{\mu_\Cbb}\)\sup_{\bsz \in \Ocal} 
C_t\(h_{l}^t \|u(\bsz)\|_{\Xcal_t} + h_{l-1}^t \|u(\bsz)\|_{\Xcal_t}\) 
= 
C_{t,\mu_\Cbb}' h_{l}^t\sup_{\bsz \in \Ocal} \|u(\bsz)\|_{\Xcal_t}.
\end{align*}
\end{proof}
\subsubsection{Nominal inf-sup conditions imply uniform inf-sup conditions}
\label{sec:NomUnifInfSup}
The preceding result, Theorem \ref{thm:CplxInfSup}, requires the validity of a
\emph{uniform discrete inf-sup condition} for the PG discretization; here, uniformity
is understood with respect to the discretization parameter $h>0$ and with respect to
the parameter sequence  $\bsz \in \Ocal$ in Theorem \ref{thm:CplxInfSup} or
with respect to $\bsy\in U$ in \eqref{eq:infsupA}, respectively.
In what follows, we assume that the two one-parameter families
of dense subspaces $\{ \bcX^h \}_{h>0}\subset \bcX$ and 
$\{ \bcY^h \}_{h>0}\subset \bcY$ are of equal, finite dimension 
$N^h = {\rm dim}(\bcX^h) = {\rm dim}(\bcY^h)$ and are
stable for the nominal operator $A_0\in \cL(\bcX,\bcY')$ in \eqref{eq:lineardependency},
i.e., the discrete inf-sup conditions hold
\begin{equation}\label{eq:DiscInfSupNom}
\inf_{0\ne w^h\in \bcX^h} \sup_{0\ne v^h\in \bcY^h}
\frac{ \langle A_0 w^h, v^h \rangle }{ \| w^h \|_{\bcX} \| v^h \|_{\bcY}} \geq \mu_0 >0 \;,
\quad 
\inf_{0\ne v^h\in \bcY^h} \sup_{0\ne w^h\in \bcX^h}
\frac{ \langle A_0 w^h, v^h \rangle }{ \| w^h \|_{\bcX} \| v^h \|_{\bcY}} \geq \mu_0 >0 \;.
\end{equation}

\begin{atheorem}\label{thm:uniinfsup}
Suppose that the parametric operators $A(\bsy)$, 
$A(\bsz)$ are affine-parametric, as in \eqref{eq:lineardependency}. 
Assume further that for $t\geq 0$ the sequences
$ \bsb_t = (b_{t,j})_{j\geq 1} $ in \eqref{eq:Bj} are small, in the sense that 
\eqref{eq:sumbjt} holds.
Then, \eqref{eq:sumbjt} with $t=0$ implies that the
discrete inf-sup conditions \eqref{eq:infsupA} 
hold uniformly with respect to $\bsy \in U$.

Moreover, if the sequence of polyradii $\rho = \(\rho_j\)_{j \geq 1}$ 
is admissible,  in the sense that \eqref{eq:admissibility} holds for $t=0$ and 
for some $\delta < 1-\kappa_0$, 
then the complex-parametric $A(\bsz)$ \changes{1) satisfies the uniform inf-sup conditions
\eqref{eq:infsupAcomplex} 
for $\bsz \in \Dcal_\rho = \bigotimes_{j \geq 1}\Dcal_{\rho_j}$, 
where $\Dcal_{\rho_j} := \{ z \in \Cbb: |z| \leq \rho_j \}$, and 2) is holomorphic with respect to any finite set of variables in $\Dcal_\rho$}.

Similarly, $A(\bsz)$ is invertible in $\Lcal(\Xcal_t,\Ycal_t')$ 
uniformly  for $\bsz \in \Dcal_{\rho}$ if 
$\rho$ is $\delta_t$-admissible w.r.t. the sequence 
$\bsb_t = \(b_{t,j} \)_{j \geq 1}$\changes{ with $\delta_t < 1-\kappa_t$}, where $b_{t,j} := \|A_0^{-1}A_j\|_{\Lcal(\Xcal_t)}$.
\end{atheorem}
%
\begin{proof}
Let $\bsb_0$ be such that 
condition~\eqref{eq:sumbjt} holds with $t=0$.
Since $A_0$ is assumed to be boundedly invertible, we can write 
$A(\bsy) = A_0\( I + \sum_{j \geq 1}y_j A_0^{-1}A_j \)$ and estimate
\begin{equation*}
\left\|\sum_{j \geq 1}y_j A_0^{-1}A_j\right\|_{\Lcal(\Xcal)} \leq \sum_{j \geq 1}|y_j|b_{0,j} \leq \sum_{j \geq 1} b_{0,j} := \kappa_0 < 1.
\end{equation*}
It follows from a perturbation (Neumann series) argument that the operator $A(\bsy)$ is 
uniformly boundedly invertible.
The discrete $\inf-\sup$ conditions hold with $\mu \leq \mu_0(1-\kappa_0)$. 

One may extend this argument to the complexified operator $A(\bsz)$ 
defined for $\bsz \in \Dcal_\rho$. 
This yields 
\begin{equation*}
\left\|\sum_{j \geq 1}z_j A_0^{-1}A_j\right\|_{\Lcal(\Xcal)} 
\leq \sum_{j \geq 1}|z_j|b_{0,j} 
\leq \sum_{j \geq 1} \rho_j b_{0,j} := \delta +\kappa_0 < 1
\;.
\end{equation*}
Therefore, the complex $\inf-\sup$ 
conditions~\eqref{eq:infsupAcomplex} hold with constant $\mu_{\Cbb} \leq \mu_0(\delta+\kappa)$.

The proof of the uniform invertibility in $\Lcal(\Xcal_t,\Ycal_t')$ follows in a similar fashion.

\changes{The operator $A(\bsz)$ being invertible, we may write, for $\bsz \in \Dcal_\rho$ and some $k \in \Nbb$, 
\begin{equation*}
u(\bsz) = \left( A_0 + \sum_{j \neq k}A_0 z_j A_0^{-1}A_j + A_0 z_k A_0^{-1}A_k \right)^{-1}f = \left( I + \sum_{j \neq k}z_j B_j + z_k B_k \right)^{-1}A_0^{-1}f. 
\end{equation*}
Whence, $u(\bsz)$ is holomorphic with respect to $z_k \in \Dcal_{\rho_k}$ as the image of $f$ via a resolvent operator. 
Hartogs' theorem concludes the holomorphy with respect to any finite set of parameters.}
\end{proof}

%
\subsubsection{Norm bounds on the \Tscheb gpc coefficients}
\label{sec:NormCoef}
We now estimate the magnitudes of the \Tscheb coefficients. 
These estimates are used in the next section to show the $\ell_{\omega,p}$ 
summability of the sequence of \Tscheb coefficients. 
We recall that $\Ecal_\rho = \bigotimes_{j \geq 1}\Ecal_{\rho_j}$ 
is a product of Bernstein ellipses  \changes{$\Ecal_{\rho_j} = \{ (z+z^{-1})/2, z \in \C : |z| =\rho_j\}$ 
and let $E_\rho = \bigotimes_{j \geq 1}E_{\rho_j}$ be the product of the open regions 
$E_{\rho_j} := \{ (z+z^{-1})/2, z \in \Cbb: 1 \leq |z| < \rho_j \}$ bounded by the Bernstein ellipses $\Ecal_{\rho_j}$. 
We note that $E_\rho$ and $\Ecal_\rho$ are contained in $\Dcal_\rho$ so that in particular under the assumptions of Theorem~\ref{thm:uniinfsup} we are in the setting of the next result.}

\begin{atheorem}
\label{thm:boundGPCcoef}
Let $\nu \in \Fcal$. 
Assume that the discretization spaces have the approximation property~\eqref{eq:approxProperty}. 
Additionaly, assume that there exists a sequence $\rho = \(\rho_j\)_{j \geq 1}$, with $\rho_j > 1$
such that the complex extension $\bsz \mapsto \dfrak u^{l}(\bsz)$ 
is holomorphic with respect to any finite set of variables on 
$E_\rho$ and with $A(\bsz) \in \Lcal(\Xcal_t,\Ycal_t')$ 
being uniformly boundedly invertible for every $\bsz \in \Ecal_\rho$. 
Then the \Tscheb coefficients of the difference $\dfrak u^{l} = u^l - u^{l-1}$ can be estimated as
\begin{equation*}
\|\dfrak u_{\nu}^{l}\|_\Xcal \leq C h_l^{t'}\|f\|_{\Ycal_t'}\rho^{-\nu}\;.
\end{equation*}
If in addition we assume smoothness for the functional, 
i.e. $\Gcal \in \Xcal_{t'}'$ for some $0 < t' \leq \bar{t}$, 
then it holds
\begin{equation}\label{norm:bound:dFnu}
|\dfrak F_{\nu}^{l}| \leq C h_l^{t+t'}\|f\|_{\Ycal_t'}\|\Gcal\|_{\Xcal_{t'}'}\rho^{-\nu},
\end{equation}
where the constants
depend on the smoothness parameters $t$ and $t'$ but not on $h_l$. 
\end{atheorem}

\begin{proof}
The proof is similar to the one in~\cite{Rauhut14CSPG} with 
appropriate modifications due to the introduction of the levels. 
The tensorized \Tscheb polynomials being orthogonal, 
it holds
\begin{equation*}
\dfrak u_\nu^{l} = \int \limits_U \dfrak u^{l}(\bsy) T_\nu(\bsy) \dif{\eta}(\bsy).
\end{equation*}
Consider the multi-index 
$\nu = n \ebf_1 = (n,0,0 \cdots) \in \Fcal$ and
split the parameter space as $U = [-1,1] \times U'$, 
then any parameter sequence $\bsy$
can be written as $\bsy = (y_1,\bsy')$ with $y_1 \in [-1,1]$. 
Thus 
\begin{equation}
\label{eq:detailine1}
\dfrak u_{n\ebf_1}^{l} 
= 
\int \limits_{U'} \int \limits_{-1}^{+1}T_n(t) 
\dfrak u^{l}(t,\bsy') \frac{\dif{t}}{\pi\sqrt{1-t^2}}\dif{\eta}(\bsy').
\end{equation}
With the change of variables $t = \cos(\phi)$ we obtain
\begin{align*}
\int\limits_{-1}^{+1}T_n(t) \dfrak u^{l}(t,\bsy')\frac{\dif{t}}{\pi\sqrt{1-t^2}} 
&= \frac{\sqrt{2}}{\pi}\int\limits_0^{\pi}\cos(n\phi) \dfrak u^{l}(\cos(\phi),\bsy')\dif{\phi} 
= \frac{1}{\sqrt{2}\pi} \int \limits_{-\pi}^{+\pi}\cos(n\phi)\dfrak u^{l}(\cos(\phi),\bsy')\dif{\phi}.
\end{align*}
This gives
\begin{align*}
\int\limits_{-1}^{+1}&T_n(t) \dfrak u^{l}(t,\bsy')\frac{\dif{t}}{\pi\sqrt{1-t^2}} 
= 
\frac{1}{\sqrt{2}\pi i} \int \limits_{|z|=1}\frac{z^n+z^{-n}}{2} \dfrak u^{l}\(\frac{z+z^{-1}}{2},{\bsy}'\)\frac{\dif{z}}{z} \nonumber \\ 
&= \frac{1}{2\sqrt{2}i\pi}\int \limits_{|z|=1}z^{n-1} \dfrak u^{l}\(\frac{z+z^{-1}}{2},\bsy'\)\dif{z}
 + 
\frac{1}{2\sqrt{2}i\pi}\int \limits_{|z|=1}z^{-n-1} \dfrak u^{l}\(\frac{z+z^{-1}}{2},\bsy'\)\dif{z}.
\end{align*}
Due to the assumption that the extension $\bsz \to \dfrak u^l(\bsz)$ to $E_\rho$
is holomorphic, the mappings
\begin{equation*}
z \mapsto z^{n-1}\dfrak u^{l}\(\frac{z+z^{-1}}{2},\bsy'\), \quad 
\text{ and } \quad 
z \mapsto z^{-n-1} \dfrak u^{l}\(\frac{z+z^{-1}}{2},\bsy'\)
\end{equation*}
are analytic on $E_{\rho_1}$. 
By Cauchy's theorem it follows, for $1 < \sigma < \rho_1$, that 
\begin{align*}
\int\limits_{-1}^{+1}T_n(t)\dfrak u^{l}(t,\bsy')\frac{\dif{t}}{\pi\sqrt{1-t^2}} 
& = \frac{1}{2\sqrt{2} i\pi}
\int \limits_{|z|=\sigma^{-1}}z^{n-1} \dfrak u^{l}\(\frac{z+z^{-1}}{2},\bsy'\)\dif{z} 
\\
& 
+ \frac{1}{2\sqrt{2} i\pi}\int \limits_{|z|=\sigma}z^{-n-1} \dfrak u^{l}\(\frac{z+z^{-1}}{2},
\bsy'\)\dif{z}.
\end{align*}
Now notice that $z \mapsto \dfrak u^{l}(z,\bsy')$ is bounded 
by $C' h_l^t\|f\|_{\Ycal_t'}$ (in $\Xcal$) in a polydisc containing in $E_{\rho}$.
Indeed, 
the approximation property of the discretization spaces, see Corollary~\ref{cor:normDiff}, 
together with the bounded invertibility in the smoothness spaces, ensures 
\begin{equation}
\label{eq:boundM}
\sup_{\bsz \in \Ecal_\rho}\|\dfrak u^{l}(\bsz)\|_{\Xcal} 
= \sup_{\bsz \in \Ecal_\rho} \| u^{l}(\bsz)-u^{l-1}(\bsz) \|_{\Xcal} 
\leq C h_l^t \sup_{\bsz \in \Ecal_\rho} \|u(\bsz)\|_{\Xcal_t} 
\leq C' h_l^t\|f\|_{\Ycal_t'}.
\end{equation}
It follows that 
\begin{align}
\left\| 
\int\limits_{-1}^{+1}T_n(t) \dfrak u^{l}(t,\bsy')\frac{\dif{t}}{\pi\sqrt{1-t^2}} \right\|_{\Xcal} 
&\leq \frac{1}{2\sqrt{2} \pi}\int \limits_{|z|=\sigma^{-1}} |z^{n-1}|
\left\| \dfrak u^{l}\(\frac{z+z^{-1}}{2},\bsy'\) \right\|_\Xcal\dif{z} 
\nonumber\\ 
&~+ 
\frac{1}{2\sqrt{2} \pi}
\int \limits_{|z|=\sigma}|z^{-n-1}| \left\|\dfrak u^{l}\(\frac{z+z^{-1}}{2},\bsy'\)
\right\|_\Xcal\dif{z} 
\nonumber \\ 
& \leq 
\frac{1}{2\sqrt{2} \pi \sigma^{n-1}} 2\pi C' h_l^t\|f\|_{\Ycal_t'} \sigma^{-1} 
+ \frac{1}{2\sqrt{2} \pi \sigma^{n+1}} 2 \pi \sigma C' h_l^t\|f\|_{\Ycal_t'} \nonumber\\
& = \sqrt{2}C' h_l^t\|f\|_{\Ycal_t'} \sigma^{-n}. 
\label{eq:errDiffLvlL}
\end{align}
This bound is valid for any $\sigma < \rho_1$ and hence holds up to $\sigma = \rho_1$. 

%
Finally, inserting Eq.~\eqref{eq:errDiffLvlL} back into
Eq.~\eqref{eq:detailine1} after integrating over 
$\bsy' \in U'$ with respect to the probability measure $\dif{\eta}(\bsy')$ yields 
\begin{equation*}
\|\dfrak u_{n\ebf_1}^{l}\|_{\Xcal} \leq C'h_l^t\|f\|_{\Ycal_t'}\rho^{-n}.
\end{equation*}
Similarly, given any $\nu \in \Fcal$, it follows that
\begin{equation*}
\|\dfrak u_{\nu}^{l}\|_{\Xcal} \leq C'h_l^t\|f\|_{\Ycal_t'}\rho^{-\nu},
\end{equation*}
by applying Cauchy's integral formula in $\IC$
with respect to each variable $z_j$ for $j\in \{ j: \nu_j \ne 0\}$.

The \Tscheb coefficients of the functional are estimated in a similar manner, using~\eqref{eq:errorSLFEM},
\begin{equation}
\label{eq:estimatesGnu}
|\dfrak F_{\nu}^{l}| \leq C'h_l^{t+t'}\|f\|_{\Ycal_t'}\|\Gcal\|_{\Xcal_{t'}'}\rho^{-\nu}.
\end{equation}
\end{proof}
\subsubsection{Summability of the sequence of \Tscheb gpc coefficients}
\label{sec:SeqSumm}
It remains to prove the existence of a $\delta$-admissible polyradii
$\rho$ (depending on $\nu$) and 
to verify the $\ell_{\omega,p}$-summability of the right hand side 
of \eqref{norm:bound:dFnu} with respect to $\nu \in \Fcal$, i.e.,
of the sequence $(\rho(\nu)^{-\nu})_{\nu \in \Fcal}$. 
Hereby, we identify suitable weights $\omega = (\omega_\nu)_{\nu \in \Fcal}$ as well.
%
%
In contrast to unweighted $\ell_p$-summability \cite{cohen2010convergence,cohen2011analytic}, 
weighted $\ell_{\omega,p}$-summability -- considered first in \cite{Rauhut14CSPG} -- 
requires stronger assumptions on the sequence $\(b_{0,j}\)_{j \in \Nbb}$ 
used as base for the $\delta$-admissibility~\eqref{eq:admissibility}.
Namely, with $\vbf = \(v_j\)_{j \in \Nbb}$ and $v_j \geq 1$, we ask 
for properties \eqref{eq:sumbjvjtilda} and \eqref{eq:psumbjvj} to be valid. 
\begin{atheorem}
\label{thm:wlpRho}
Let $\vbf$ be a sequence of weights fulfilling the summability 
conditions~\eqref{eq:sumbjvjtilda} and~\eqref{eq:psumbjvj} 
and let $\omega_\nu := \theta^{\|\nu\|_{0}}\vbf^{\nu}$ for any $\nu \in \Fcal$ 
and some $\theta \geq 1$. 
There exists a sequence of polyradii ${\(\bf \rho(\nu)\)_{\nu \in \Fcal}}$ 
such that 
\begin{itemize}
\item[i)] 
for each $\nu\in \Fcal$, 
${\bf \rho}= {\bf \rho}(\nu)$ is $\delta$-admissible, 
with $\delta = (1-\kappa_{v,p})/2$,  and 
\item[ii)] 
$\|\({\bf \rho(\nu)}^{-\nu}\)_{\nu \in \Fcal}\|_{\omega,p} \leq K_{\theta,p} < \infty$.
\end{itemize}
\end{atheorem}

\begin{proof}
Full details of the argument can be found in~\cite{Rauhut14CSPG}; 
here, we only indicate the main steps,
in particular the construction of a sequence of weights $\omega$ 
and an associated, admissible sequence of polyradii.

For the weights $\vbf$ and a constant $\theta \geq 1$, we define the sequence of weights
\begin{equation}
\label{eq:omegas}
\omega_\nu(\theta) := \theta^{\|\nu\|_0}\vbf^\nu = \theta^{\|\nu\|_0}\prod_{j: \nu_j \neq 0}v_j^{\nu_j},\quad \text{for all } \nu \in \Fcal.
\end{equation}
Because of \eqref{eq:sumbjvjtilda}, there exists a finite set $E \subset \Nbb$ such that, with $F := \Nbb \backslash E$, 
\begin{equation*}
\sum_{j \in F} v_j^{(2-p)/p}b_{0,j} \leq \frac{\delta}{8 \theta^{(2-p)/p}}.
\end{equation*}
For a given constant $\alpha > 1$ with $(\alpha - 1)\sum_{j \in E}v_j^{(2-p)/p}b_{0,j} < \delta/2$, we define the sequences of polyradii (generally depending on $\nu$) as 
\begin{equation}
\label{eq:rhoofnu}
\rho_j(\nu) = \left \{ \begin{array}{cl} \alpha v_j^{(2-p)/p}, &j \in E, \\ \operatorname{max}\left\{ v_j^{(2-p)/p}, \frac{\nu_j}{2|\nu_F| b_j} \right\}, &j \in F, \end{array} \right.
\end{equation}
where we used the notation $|\nu_F| := \sum_{j \in F}\nu_j$.
The $\delta$-admissibility of this sequence, as well as its $\ell_{\omega,p}$ summability,
ensuring the summability of the \Tscheb expansion of the differences, 
have been proved in~\cite[Theorem 4.2]{Rauhut14CSPG}.
\end{proof}
Combining the estimate~\eqref{eq:estimatesGnu} with the 
$\ell_{\omega,p}$ summability of the sequence $\rho$ yields
\begin{equation}
 \label{eq:estimateNormG}
 \|\dfrak F^{l}\|_{\omega,p} \leq Ch_l^{t+t'}\|f\|_{\Ycal_t'}\|\Gcal\|_{\Xcal_{t'}'}\|\(\rho^{-\nu}\)_{\nu \in \Fcal}\|_{\omega,p}\;.
\end{equation}
Consequently, with Eq.~\eqref{eq:approxSampPtsTscheb} it follows
\begin{equation}
\label{eq:linfLvlldfrak}
 |\dfrak F^{l}(\bsy) - \ucs{\deltau{F}^l}(\bsy)| 
\leq 
C s_l^{1-1/p} h_l^{t+t'}\|f\|_{\Ycal_t'}\|\Gcal\|_{\Xcal_{t'}'}\|\(\rho(\nu)^{-\nu}\)_{\nu \in \Fcal}\|_{\omega,p}
\;.
\end{equation}
Theorem~\ref{thm:summability} is a direct consequence of the results in this section. 
Indeed, the bounded invertibility in the smoothness spaces 
of $A_0 \in \Lcal(\Xcal_t,\Ycal_t')$ together with the 
summability~\eqref{eq:deltatAdmissibility} implies the uniform bounded 
invertibility of the operator $A(\bsz) \in \Lcal(\Xcal_t,\Ycal_t')$, 
via a perturbation argument as stated in Theorem~\ref{thm:uniinfsup}. 
This ensures the applicability of Theorem~\ref{thm:boundGPCcoef} 
(which itself depends on the two previous theorems). 
Theorem~\ref{thm:wlpRho} finally proves the existence of both a positive
weight sequence $\omega$ and a sequence of polyradii $\rho$  
as well as the $\ell_{\omega,p}$ summability. 
\subsection{Rate of convergence of the MLCSPG method}
\label{sec:convergence}
To simplify the exposition, 
we only derive the bounds for the approximation of a functional of the parametric solution. 
The results can be applied mutatis mutandis to derive the 
convergence rates for the full solution $u(\bsy)$, 
-- once the details of the (single level) compressive sensing
scheme for the approximation of the full solution are worked out.
We continue the estimate in~\eqref{eq:errorToAnalyze} as follows:
\begin{align*}
|F(\bsy)-F^L_{\text{MLCS}}(\bsy)| 
&\leq |F(\bsy) - F^L(\bsy)| + \sum_{l=1}^L \left| \dfrak F^{l}(\bsy) - \ucs{\deltau{F}^l}(\bsy)\right| \\ 
&\leq Ch_{L}^{t+t'}\| f \|_{\Ycal_t'}\| \Gcal \|_{\Xcal_{t'}'} 
+ \sum_{l=1}^L C s_l^{1-1/p}\vertiii{\dfrak F^{l}}_{\omega,p} 
\\ 
&\leq C\| f \|_{\Ycal_t'}\| \Gcal \|_{\Xcal_{t'}'} 
\( h_{L}^{t+t'} + \sum_{l=1}^L s_l^{1-1/p}h_{l}^{t+t'} \|\(\rho(\nu)^{-\nu}\)_{\nu \in \Fcal}\|_{\omega,p}\).
\end{align*}
We absorb the norm $\|\(\rho(\nu)^{-\nu}\)_{\nu \in \Fcal}\|_{\omega,p}$ into the constant $C>0$, yielding
\begin{equation*}
|F(\bsy)-F^L_{\text{MLCS}}(\bsy)| 
\leq 
C\| f \|_{\Ycal_t'}
\| \Gcal \|_{\Xcal_{t'}'} 
\( h_{L}^{t+t'} + \sum_{l=1}^L s_l^{1-1/p}h_{l}^{t+t'}\)
\;.
\end{equation*}

Using that the levels are related via 
$h_l = h_{l-1}/2$ we obtain
\begin{equation*}
|F(\bsy)-F^L_{\text{MLCS}}(\bsy)| 
\leq 
C\| f \|_{\Ycal_t'}\| \Gcal \|_{\Xcal_{t'}'} h_{L}^{t+t'} \( 1 + \sum_{l=1}^L s_l^{1-1/p}2^{(L-l)(t+t')}\)
\;.
\end{equation*}
We balance sampling and discretization errors on each mesh level $l$ in this bound.
Thus the choice
\begin{equation} \label{eq:sl} 
s_l \asymp 2^{(L-l)(t+t')p/(1-p)} 
= 
2^{(L-l) \sigma_p(t+t')}, \qquad \mbox{ with } \sigma_p(t) = \frac{tp}{1-p},
\end{equation}
implies an overall error bound of 
\begin{equation*}
|F(\bsy)-F^L_{\text{MLCS}}(\bsy)| 
\leq 
C\| f \|_{\Ycal_t'}\| \Gcal \|_{\Xcal_{t'}'} (L+1) h_{L}^{t+t'} = C\| f \|_{\Ycal_t'}\| \Gcal \|_{\Xcal_{t'}'} (|\log(h_L)|+1) h_{L}^{t+t'}
\;.
\end{equation*}
\changes{From the choice of the sparsities~\eqref{eq:sl} 
together with \eqref{eq:mwl1approx}, it follows that a number of samples per level scaling as}
\begin{equation}
\label{eq:ml}
m_l \asymp s_l \max\{ \log^3(s_l) \log(N_l), \log(1/\gamma)\} \asymp 2^{(L-l)\sigma_p(t+t')}\max\{ \(L-l\)^3 \log(N_l), \log(1/\gamma)\}
\end{equation}
is sufficient for the error bound~\eqref{eq:linfLvlldfrak} to be valid at level $l$ with probability exceeding \changes{$1 - \gamma$}.
Note that the size $N_l$ of the initial index set $\Gamma_l$ 
may depend on $s_l$ and on the choice of weights $\omega$. 
More details are given in the next section.
The global error in $L^2$ is bounded as in Eq.~\eqref{eq:errorToAnalyze}, 
\begin{align}
\label{eq:errorAnalyzeL2}
\|F-F^L_{\text{MLCS}}\|_2 
\leq 
\|F-{F^L}\|_2 + \sum_{l = 1}^L \|\dfrak F^l - \ucs{\deltau{F^l}}\|_2
\;.
\end{align}
The first term is computed using the 
uniform bound~\eqref{eq:errorSLFEM} and the fact that $\eta$ is a probability measure. 
To compute the sum, it suffices to apply 
the $\ell_2$-error bound in \eqref{eq:approxSampPts} to the $L$ details, 
$\| \dfrak F^l - \ucs{\deltau{F^l}} \|_2 \leq D s_l^{1/2-1/p}\|\dfrak F^l\|_{\omega,p}$.
Hence, applying \eqref{eq:errorSLFEM} to the first term and combining the $\ell_2$ bound in 
\eqref{eq:approxSampPts} with the prior estimate~\eqref{eq:priorBoundFnu} 
and the number of samples~\eqref{eq:sl} in the terms in the sum yields
$$
\|F-F^L_{\text{MLCS}}\|_2 
\leq 
C\|f\|_{\Ycal_t'}\|\Gcal\|_{\Xcal_{t'}'}
h_L^{t+t'}\(1 + \sum_{l=1}^L 2^{(t+t')(L-l)\frac{p-2}{2(1-p)}}2^{(L-l)(t+t')}\) 
=
C' 
\|f\|_{\Ycal_t'}\|\Gcal\|_{\Xcal_{t'}'}h_L^{t+t'}
\;.
$$
Alternatively, one can also balance the number of samples with 
the discretization error to reach a prescribed $L_2$ error of 
$\cO(h_L^{t+t'})$ by combining Eq.~\eqref{eq:errorAnalyzeL2} 
with the compressed sensing approximation~\eqref{eq:approxSampPts}: 
\begin{align*}
\|F-F^L_{\text{MLCS}}\|_2 
&\leq Ch_L^{t+t'}\|f\|_{\Ycal_t'}\|\Gcal\|_{\Xcal_{t'}'} 
       + C \sum_{l = 1}^L s_l^{1/2-1/p}h_l^{t+t'}\|f\|_{\Ycal_t'}\|\Gcal\|_{\Xcal_{t'}'}\|\(\rho(\nu)^{-\nu}\)_{\nu \in \Fcal}\|_{\omega,p} 
\nonumber \\
&\leq C \|f\|_{\Ycal_t'}\|\Gcal\|_{\Xcal_{t'}'} h_L^{t+t'} \( 1 + \sum_{l=1}^L s_l^{1/2-1/p} 2^{(L-l)(t+t')} \)
\;.
\end{align*}
In this case, choosing
\begin{equation}
\label{eq:sl2}
s_l \asymp 2^{\frac{(L-l)(t+t')2p}{2-p}}
\end{equation}
ensures the $L_2$ error bound 
\begin{equation}
\label{eq:MLErBd}
\|F-F^L_{\text{MLCS}}\|_2 
\leq C 
h_L^{t+t'} (1+|\log h_L|) \|f\|_{\Ycal_t'}\|\Gcal\|_{\Xcal_{t'}'} 
\;.
\end{equation}
\section{Implementation Aspects}
\label{sec:implem}
This section describes several aspects that are relevant 
for the numerical applicability of the theoretical approach introduced above.
In particular, we investigate the truncation of the (potentially infinite) 
sequence of parameters to a finite subset,
and specify initial choices of finite index sets $\Lambda \subset \Fcal$
that are guaranteed to contain the support of the best (weighted) $s$-term approximation 
of the solution and can be used within weighted $\ell_1$-minimization or other CS algorithms.
%
\subsection{Dimension truncation}
\label{ssec:truncation}
So far, we have worked on a purely theoretical level, where the parameter vector is potentially infinite (but countable). 
To ensure the applicability of the results, we have to verify that truncating the parameter vector to a finite dimensional space (yet allowing this truncation to be rather large) still allows for reliable approximations. 

We consider the weak solutions of the truncated version of Eq.~\eqref{eq:operator}:
\begin{equation}
\label{eq:weakTruncated}
\text{Find $u^{(B)} \in \Xcal$, such that } \langle A^{(B)}(\bsy)u^{(B)}, v \rangle = \langle f, v \rangle \quad \mbox{ for all } v \in \Ycal,
\end{equation}
where the operator $A^{(B)}(\bsy)$ is defined, for a finite $B \in \Nbb$, as $A(y_1, y_2, \cdots, y_B, 0, 0, \cdots)$.

In particular, we assume some decay of the \emph{energy} of the operator $A(\bsy)$ 
(i.e. assuming a certain order on the parameters) such that for any $\varepsilon > 0$, 
there exists $B := B(\varepsilon,A)$ with 
\begin{equation}
\label{eq:decayEnergyOperator}
\|A(\bsy) - A^{(B)}(\bsy)\|_{\Lcal(\Xcal, \Ycal')} \leq \varepsilon \mu, \quad \forall \bsy \in U,
\end{equation}
where $\mu$ is the constant appearing in the $\inf-\sup$ conditions~\eqref{eq:infsupA}. 

In this case, the following generalization of results in~\cite{Dick13QMCPG} holds.
\begin{aprop}
\label{prop:truncation}
Assume the operator $A$ satisfies the (continuous) $\inf-\sup$ conditions~\eqref{eq:infsupA} 
and the decay property~\eqref{eq:decayEnergyOperator}. 
Then for any accuracy parameter $\varepsilon$, 
there exists a truncation parameter $B \in \Nbb$ such that the solutions 
to the truncated problem~\eqref{eq:weakTruncated} and 
to the original problem~\eqref{eq:operator} are close to each other in the following sense
\begin{align}
\label{eq:approxDecayEnergyOperator}
\|u^{(B)}(\bsy) - u(\bsy)\|_{\Xcal} &\leq \frac{C \varepsilon}{\mu}\|f\|_{\Ycal'}, 
\end{align}
where $u^{(B)}(\bsy)$ is the solution of the truncated problem~\eqref{eq:weakTruncated}.
\end{aprop}

\begin{proof}
The weak solutions are characterized by
\begin{align*}
\text{Find $u(\bsy)$ such that } \quad \langle A(\bsy)u(\bsy), v\rangle &= \langle f, v\rangle \quad \mbox{ for all } v \in \Ycal, \\
\text{Find $u^{(B)}(\bsy)$ such that } \quad  \langle A^{(B)}(\bsy)u^{(B)}(\bsy), v\rangle &= \langle f, v\rangle \quad \mbox{ for all } v \in \Ycal.
\end{align*}
Since these equalities hold for all $v$, they imply the orthogonality conditions
\begin{equation*}
\langle A(\bsy)u(\bsy) - A^{(B)}(\bsy)u^{(B)}(\bsy), v\rangle = 0 \quad  \mbox{ for all } v \in \Ycal. 
\end{equation*}
Rearranging the terms yields 
$\langle A(\bsy)\(u(\bsy) - u^{(B)}(\bsy)\), v \rangle 
  = 
  -\langle\(A(\bsy)-A^{(B)}(\bsy)\)u^{(B)}(\bsy),v\rangle$ for all $v \in \Ycal$. 
This means that $u(\bsy) - u^{(B)}(\bsy)$ is the weak solution 
to the operator equation~\eqref{eq:operator} with forcing term 
$\(A(\bsy)-A^{(B)}(\bsy)\)u^{(B)}(\bsy)$. 
Consequently, using the $\inf-\sup$ conditions twice 
and the decay property~\eqref{eq:decayEnergyOperator}, we obtain
\begin{equation*}
\|u(\bsy) - u^{(B)}(\bsy)\|_{\Xcal} 
\leq 
\frac{C}{\mu}\left\|A(\bsy)-A^{(B)}(\bsy)\right\|_{\Lcal(\Xcal,\Ycal')} \|u^{(B)}(\bsy)\|_{\Xcal} 
\leq \frac{\varepsilon C}{\mu}\|f\|_{\Ycal'},
\end{equation*}
which concludes the proof.
\end{proof}
Consequently, 
it is sufficient to draw the $m_l$ samples per level at random according to the truncated distribution. 
As a concrete example let us consider the case of linear dependence
on the parameters as described in~\cite{Rauhut14CSPG} and 
in Eq.~\eqref{eq:lineardependency}. 
Assuming that $A_0 : \Xcal \to \Ycal'$ is invertible 
(which was required in Theorem~\ref{thm:summability}) and that 
$(b_{0,j})_j \in \ell_1$ (which is weaker than the conditions in the previous section) 
the fluctuations $A_j$, $j \geq 1$ are arranged in nonincreasing order, i.e., 
such that $b_{0,j} \geq b_{0,k}$ for $1 \leq j \leq k$, 
then the 
operator~\eqref{eq:lineardependency} satisfies the following
dimension truncation error bound
\begin{equation*}
\|A(\bsy) - A^{(B)}(\bsy)\|_{\Lcal(\Xcal,\Ycal')} 
= 
\|\sum_{j > B}y_jA_j\|_{\Lcal(\Xcal,\Ycal')} 
= 
\|A_0 \sum_{j > B}y_jB_j\|_{\Lcal(\Xcal,\Ycal')} 
\leq 
\|A_0\|_{\Lcal(\Xcal,\Ycal')}\sum_{j > B}b_{0,j},
\end{equation*}
for any $\bsy \in U$.
Moreover (see~\cite[Thm 2.9]{Rauhut14CSPG},\cite[Thm 5.1]{Kuo2012qmcRdmPDE}), the tail 
can be estimated by 
\begin{equation*}
\sum_{j > B}b_{0,j} \leq \min\left\{\frac{1}{1/p-1}, 1\right\}\|(b_{0,j})_j\|_p B^{-(1/p-1)}
\end{equation*}
for some $p < 1$. Consequently, choosing $B \geq h_L^{-(t+t')p/(1-p)}$ yields a global approximation (accounting for the truncation error, the PG approximation error, and the CS error) in $\Ocal(h_L^{t+t'})$.
%
\subsection{Initial set of candidate vectors}

As detailed in the discussion before Theorem~\ref{thm:summability}, 
the results are, so far, developed for an infinite \Tscheb expansion. 
To render the problem computationally feasible, we
truncate to a finite-dimensional, parametric expansion, where the 
truncation dimension is at our disposal and therefore can be considered
a discretization parameter.
Let the sums~\eqref{eq:expUdl} and~\eqref{eq:expFdl} be truncated to a finite set $\Gamma_l \subset \Fcal$. 
Some strategies for selecting such a set $\Gamma_l$ were already described 
in~\cite{Rauhut14CSPG}, which was based on the work in~\cite{Rauhut13wCS}.
We have the following analog to Theorem~\ref{thm:wl1approx} (proven in~\cite{Rauhut13wCS})
in the case of expansions in terms of a countable sequence of parameters.
\begin{atheorem}
\label{thm:wl1Truncated}
Let $\gamma \in (0,1)$. 
Let $F(\bsy) = \sum_{\nu \in \Fcal}F_\nu T_\nu(\bsy)$ 
be a function with
$\|{\bf F}\|_{\omega,p} = \left\|\( F_\nu\)_{\nu}\right\|_{\omega,p} < \infty$ 
for some $p < 1$ and some weights $\omega_\nu \geq \|T_\nu\|_\infty$ for all $\nu \in \Fcal$. 
For a given sparsity $s_l \geq 1$, define the initial set of indices 
as
\begin{equation}
 \label{eq:gamma0}
 \Gamma_l := \{\nu \in \Fcal: \omega_\nu^2 \leq s_l/2 \}.
\end{equation}
Furthermore, assume that $N_l := |\Gamma_l|$ is finite and draw 
\begin{equation}\label{bound:ml}
 m_l \geq c_0 s_l \max\{\log^3(s_l)\log(N_l), \log(1/\gamma)\}
\end{equation}
sampling points $\bsy^{(i)}$ 
independently and identically distributed according to the orthogonalization measure $\eta$. 
Let $\ucs{{\bf F}}$ 
be the solution of 
\begin{equation*}
 \min \| \mathbf{H} \|_{\omega,1} \quad \text{subject to } \|A\mathbf{H} - b\|_2 \leq 2^{1-p}\tau \sqrt{m_l}s_l^{1/2-1/p}\|{\bf F}\|_{\omega,p},
\end{equation*}
for some $\tau \geq 1$ and set $\ucs{F} = \sum_{\nu \in \Gamma_l} \ucs{{\bf F}}_\bnu T_\bnu$. 
Then, with probability at least $1-\gamma$ 
\begin{align*}
 \|F - \ucs{F}\|_\infty &\leq \| {\bf F} - \ucs{{\bf F}}\|_{\omega,1} 
 \leq c_\tau s_l^{1-1/p}\|{\bf F}\|_{\omega,p}, 
\\
\|F - \ucs{F}\|_2 
&= \| {\bf F} - \ucs{{\bf F}}\|_{2} \leq d_\tau s_l^{1/2-1/p}\|{\bf F}\|_{\omega,p}
\;.
\end{align*}
\end{atheorem}
A drawback of the recovery based on an optimization problem 
is that it requires the knowledge (or an approximation) of the norm of the unknown vector ${\bf F}$. 
This can be overcome in practical applications by applying the recovery to various estimations (similar to a cross validation in the machine learning literature~\cite{ward2009CSCV}) or by using 
greedy methods, e.g.~\cite{Bykowski15MSc,Fell15wIHT}. 

The cardinality $N_l$ of the set 
\[
\Gamma_l = \{\nu \in \Fcal: \omega_\nu^2 \leq s_l/2 \} = \{ \nu \in \Fcal  : \|\nu\|_0 \log(\theta) + \sum_{j \in \operatorname{supp} \nu} 2\log(v_j) \nu_j \leq \log(s_l/2)  \},
\]
where the weights $\omega_\nu$ are chosen as in \eqref{def:weight:omega},
influences the number $m_l$ of 
samples in \eqref{bound:ml} (and the computational
complexity of the weighted $\ell_1$-minimization problem). Obviously, $N_l$ depends on $s_l$ as well as on 
the weight sequence $(v_j)$ used in the definition \eqref{def:weight:omega} of $(\omega_\nu)$. 
We recall the following estimates from \cite{Rauhut14CSPG}.

\begin{aprop}
Let $\omega_\nu = \theta^{\|\nu\|} \mathbf{v}^{\nu}$, $\nu \in \Fcal$, 
for a sequence $\mathbf{v} = (v_j)_{j \geq 1}$ specified below 
and assume $s_l \geq 1$.
\begin{enumerate}
\item 
For $v_j = \beta$ for $1 \leq j \leq d$ and $v_j = \infty$ for $j > d$ 
(i.e., we consider constant weights for the first $d$ dimensions 
and ignore the remaining ones), we have
\begin{align}
N_l = |\Gamma_l| 
&\leq \left\{ \begin{array}{ll} 
\(\(1+\frac{1}{\log_2(\beta^2)}\)ed\)^{\log_{2\beta^2}(s_l/2)}, 
& s_l < 2^{d+1}\beta^{2d}, 
\\ 
(\log_{\beta^2}(\beta^2s_l/2))^d, & s_l \geq 2^{d+1}\beta^{2d}\;. 
\end{array} \right. \nonumber 
  \end{align}
\item 
For polynomially growing weights $v_j = cj^\alpha$ with $c > 1$ and $\alpha > 0$, 
there holds subexponential growth
  \begin{align}
	 N_l & \leq C_{\alpha,c}s_l^{\gamma_{\alpha,c}\log(s_l)} \nonumber 
\end{align}
for some constants $C_{\alpha,c} > 0$ and $\gamma_{\alpha,c} > 0$ 
depending only on $c$ and $\alpha$.
 \end{enumerate}
\end{aprop}
Inserting these bounds into Condition 
\eqref{bound:ml} on the number of required samples 
(assuming that the $\log(1/\gamma)$-term does not exceed the other logarithmic terms) 
shows that the following choices of $m_l$ are valid:
\begin{itemize}
\item For constant weights $v_j = \beta$ for $1 \leq j \leq d$ and $v_j = \infty$ 
for $j > d$, we can chose
\begin{equation}
\label{eq:ml_constantv}
	m_l \asymp\left\{ \begin{array}{ll} \log(d)s_l \log^4(s_l), & s_l < 2^{d+1}\beta^{2d}, \\ ds_l\log^3(s_l)\log(\log(s_l)), & s_l \geq 2^{d+1}\beta^{2d}. \end{array} \right.
\end{equation}
\item For polynomially growing weights $v_j = cj^\alpha$ with $c > 1$ and $\alpha > 0$, 
we can chose
\begin{equation}
\label{eq:ml_polynomv} m_l \asymp s_l\log^5(s_l)\;.
\end{equation}
\end{itemize}

The case of exponentially growing weights has been analyzed in~\cite{Rauhut14CSPG} 
and yields situations where $N_l \leq m_l$. 
In this situation, compressed sensing techniques should not be used,
as least-squares methods are expected to perform better~\cite{Migliorati14discreteL2}.

We note that in the case of constant weights, 
the first case in \eqref{eq:ml_constantv} is the most relevant.
In fact, with the choice of $s_l$ as in~\eqref{eq:sl}, 
i.e., $s_l = \Ccal 2^{(L-l)(t+t')p/(1-p)}$ for 
some proportionality constant $\Ccal>0$, if 
$c := \frac{d+1+2d\log_{2}(\beta)}{(t+t')p}(1-p) - \frac{\log_2(\Ccal)(1-p)}{(t+t')p}$ 
is large enough (for instance $c \geq L$, which is true whenever $\Ccal \leq s^{d+1}\beta^{2d}/2^{L(t+t')p/(1-p)}$) 
then only the first case of \eqref{eq:ml_constantv} will occur for all $l=1,\hdots, L$.
In particular, with all the parameters ($\beta$, $t$, $t'$, and $p$) fixed, 
a larger number $d$ of active variables will lead to a larger $c$. 
It is therefore reasonable to assume that this corresponds to the main regime.

\subsection{Computational Cost}
\label{sec:comp:costs}
In the ensuing work bounds, we assume at our disposal \changes{multigrid}
solvers as described, e.g. in~\cite{HackbMGM,XuMGMSurvey}. 
These solvers compute approximate solutions of the Galerkin equations
at cost
scaling linearly in the number of unknowns of the mesh. 
This gives rise to the following complexity estimates, 
where we treat the case of constant and polynomially growing weights. 
\begin{aprop} 
\label{prop:computations}
Under the assumptions \eqref{eq:gamma0}, \eqref{bound:ml} as well as \eqref{eq:sumbjvjtilda}, 
\eqref{eq:psumbjvj} for some $0 < p < 1$ and smoothness parameters $t,t'$,
the function $\bsy \mapsto \Gcal(u(\bsy))$ can be approximated in $L^2(U,\eta)$ 
to accuracy $\mathcal{O}(h_L^{t+t'})$ via 
a MLCSPG discretization with $L$ levels and with total work $W_L^T$
scaling as 
\begin{align}
W_L^T \lesssim 
\left\{ \begin{array}{lll} 
\dfrac{\log(d) \sigma_p(\tau)^4 L^42^{L\sigma_p(\tau)}}{\sigma_p(\tau)-n}, 
&\sigma_p(\tau) > n, &(v_j) 
\text{ constant} 
\\
  \log(d) \sigma_p(\tau)^4 L^5 2^{nL}, & \sigma_p(\tau) = n, & (v_j) \text{ constant}
\\
  \dfrac{\log(d)\sigma_p(\tau)^4 (2^{nL} - 2^{\sigma_p(\tau)L})}{\(n-\sigma_p(\tau)\)^4}, 
    &\sigma_p(\tau) < n, & (v_j) \text{ constant} 
\\
    \dfrac{\sigma_p(\tau)^5 L^5 2^{L\sigma_p(\tau)}}{\sigma_p(\tau)-n}, &\sigma_p(\tau) > n, 
    & (v_j) \text{ polynomial} 
\\	
 \sigma_p(\tau)^5 L^6 2^{nL}, & \sigma_p(\tau) = n, & (v_j) \text{ polynomial}
\\
\dfrac{\sigma_p(\tau)^5 (2^{nL}-2^{\sigma_p(\tau)L})}{(n-\sigma_p(\tau))^5}, 
    &\sigma_p(\tau) < n, & (v_j) \text{ polynomial}
\end{array} \right.
\end{align}
where $\sigma_p(\tau) = \tau p/(1-p)$ with $\tau = t+t'$ 
and where $n$ denotes the spatial dimension.
\end{aprop}

\begin{proof}
Multigrid solvers have a computational complexity scaling linearly with the number 
$w_l \asymp 2^{nl}$ of unknowns at level $l$ which implies that the work at level $l$ is 
on the order of $W_l = m_l \cdot w_l$, $1 \leq l \leq L$.
 
Assuming we are given constant weights $v_j = \beta$, 
for $1 \leq j \leq d$, and $s_l < 2^{d+1}\beta^{2d}$, 
and that $d$ is sufficiently large, 
we can chose $m_l$ as in the first row of Eq.~\eqref{eq:ml_constantv}. 
Thus, omitting constants,
\begin{align}
W_L^T &= \sum_{l = 1}^L W_l 
\lesssim \sum_{l=1}^L \log(d)s_l\log^4(s_l)2^{nl} 
\lesssim \sum_{l = 1}^L\log(d)2^{(L-l)\sigma_p(\tau)}\log^4\(2^{(L-l)\sigma_p(\tau)}\)2^{nl} 
\nonumber \\
\label{eq:WlEstimateSmallDim}			
&\lesssim \log(d)\sigma_p(\tau)^4 2^{nL}\sum_{l = 1}^L(L-l)^4 2^{(L-l)(\sigma_p(\tau) -n)} 
= 
\log(d)\sigma_p(\tau)^4 2^{nL}\sum_{j = 1}^{L-1}j^4 2^{j(\sigma_p(\tau) -n)}.
\end{align}
We can bound 
$S := \sum_{j = 1}^{L-1}j^4 2^{j(\sigma_p(\tau) -n)} 
\leq \int_0^L2^{x(\sigma_p(\tau) - n)}x^4 \dif{x}$. 
If $\sigma_p(\tau) = n$, it follows that $S \leq L^5/5$. 
Otherwise, with $K = (\sigma_p(\tau)-n)\ln(2)$, an integration by part yields
\begin{equation}
S \leq \frac{L^4e^{LK}}{K} - \frac{4}{K}\int_0^L x^3 e^{xK}\dif{x}.
\end{equation}
If $K > 0$, i.e. $\sigma_p(\tau) > n$, the remaining integral is positive and thus 
$S \leq \frac{L^4e^{LK}}{K} = \frac{L^42^{L(\sigma_p(\tau)-n)}}{\ln(2)(\sigma_p(\tau)-n)}$. 
If $K < 0$, repeated integration by parts leads to 
\begin{equation}
S \leq \frac{L^4e^{LK}}{K}  -\frac{4L^3e^{LK}}{K^2} + \frac{12L^2e^{LK}}{K^3} - \frac{24Le^{LK}}{K^4} + \frac{24}{K^4}\int_0^L e^{xK}\dif{x}.
\end{equation}
Noticing that $\frac{L^4e^{LK}}{K}  -\frac{4L^3e^{LK}}{K^2} + \frac{12L^2e^{LK}}{K^3} - \frac{24Le^{LK}}{K^4} < 0$, it follows that 
\begin{equation}
S \leq \frac{24}{K^4}\int_0^L e^{xK}\dif{x} = 24 \frac{e^{LK}-1}{K^5} = \frac{24({1- 2^{(\sigma_p(\tau)-n)L}})}{(n-\sigma_p(\tau))^5\ln(2)^5}.
\end{equation}
The result for polynomially growing weight sequences $(v_j)$ is shown
in a similar fashion (with appropriate changes in exponents). 
\end{proof}

\begin{rmk} 
\label{remk:work}
Recalling
that the workload for the computation of one solution at the finest 
discretization level $L$ is $w_L \asymp 2^{nL}$,
the previous result means that
for $\sigma_p(t+t') < n$, the total work is bounded only by a 
multiple of the cost of one PDE solve at the finest level, 
where the multiplicative constant involves a factor of $\log(d)$ in the case
of constant weights and in addition only depends on $n,p, t, t'$.
\end{rmk}

Combining Theorem~\ref{thm:summability} together with 
Proposition~\ref{prop:computations} about the computation 
costs and Proposition~\ref{prop:truncation} regarding the truncation of the operator, 
we are finally able to state our main theorem.
To this end we first summarize the assumptions on the parametric operator 
$A(\bsy) = A_0 + \sum_{j \geq 1} y_j A_j$. 
\begin{itemize}
\item The nominal operator $A_0$ is inf-sup stable, i.e.,
\[
\inf_{0\ne w^h\in \bcX^h} \sup_{0\ne v^h\in \bcY^h}
\frac{ \langle A_0 w^h, v^h \rangle }{ \| w^h \|_{\bcX} \| v^h \|_{\bcY}} \geq \mu_0 >0,
\quad 
\inf_{0\ne v^h\in \bcY^h} \sup_{0\ne w^h\in \bcX^h}
\frac{ \langle A_0 w^h, v^h \rangle }{ \| w^h \|_{\bcX} \| v^h \|_{\bcY}} \geq \mu_0 >0.
\]
\item 
For some $0 < p < 1$ and some weight sequence $\mathbf{v} = (v_j)_{j \in \Nbb}$ with $v_j \geq 1$, 
the sequence $\bsb_0$ with components $b_{0,j} = \|A_0^{-1} A_j\|_{\mathcal{L}(\Xcal)}$, $j \geq 1$, 
satisfies
\[
\kappa_{\mathbf{v},p} 
:= 
\sum_{j \geq 1} 
b_{0,j} v_j^{(2-p)/p} < 1 \quad \mbox{ and } \sum_{j \geq 1} b_{0,j}^p v_j^{2-p} < \infty. 
\]
\item 
For some $t \in (0,\bar{t})$, the operators $A_j$, $j \geq 0$, 
are defined as operators from $\Xcal_t$ into $\Ycal_t'$  
the sequence $\mathbf{b}_t$ with components $b_{t,j} = \|A_0^{-1} A_j \|_{\mathcal{L}(\Xcal_t)}$ 
satisfies
\[
\kappa_{t} := \sum_{j \geq 1} b_{t,j} \leq 1\;,\quad 
\bsb_t \in \ell^{p_t}.
\]
\end{itemize}
\begin{atheorem}
\label{thm:mainThm}
\changes{Let $L \in \Nbb$ be a number of discretization levels and $\gamma_\ell \in (0,1)$, $\ell=1,\hdots,L$.}
Let $A(\bsy)$ be an affine-parametric operator and 
let $\vbf = (v_j)_{j \geq \Nbb}$ be a sequence of weights with $v_j \geq 1$. 
Assume that 
$A_0 \in \Lcal(\Xcal_t,\Ycal_t')$ is boundedly invertible and that 
the sequence
$\bsb_t = (b_{t,j})_{j \geq 1}$ are
such that the summability conditions~\eqref{eq:sumbjt} 
and~\eqref{eq:psumbjvj} hold true for some $0 < p < 1$. 
Then, for any discretization level $1 \leq l \leq L$, 
the sequence of \Tscheb coefficients of $\dfrak u^{l}$ with respect to 
the parameter vector~\eqref{eq:expUdl} is (weighted) compressible, 
i.e., for a sequence of weights 
$\omega = (\omega_\nu)_{\nu \in \Fcal}$ with $\omega_\nu = \theta^{\|\nu\|_0}\vbf^\nu$ 
there holds 
$\sum_{\nu \in \Fcal}\omega^{2-p} \|\dfrak u_\nu^l\|_\Xcal^p < \infty$. 

Moreover, if we are interested in a functional of the solution 
$F(\bsy) = \Gcal(u(\bsy))$ and if the operators 
$A(\bsy)\in \Lcal(\Xcal_t,\Ycal_t')$
are boundedly invertible in the smoothness scales $(\Xcal_t , \Ycal_t )$ 
in \eqref{eq:SmoothPrimal}, \eqref{eq:SmoothDual} and if
$\Gcal \in \Xcal_{t'}'$, then the function 
$F(\bsy) = \sum_{\nu \in \Fcal}F_{\nu}T_\nu(\bsy)$ can be approximated 
by $F^L_{\text{MLCS}}(\bsy) := \sum_{l=1}^L \ucs{\deltau{F}^l}(\bsy)$ 
where $\ucs{\deltau{F}^l}(\bsy)$ is a single-level CSPG approximation from 
\[
\changes{m_l \asymp s_l \max\{\log^3(s_l) \log(N_l), \log(1/\gamma_l)\}}
\]
sampling points with $s_l \asymp 2^{(L-l)(t+t')p/(1-p)}$, 
where
$N_l = |\Gamma_l|$ for $\Gamma_l = \{ \nu \in \Fcal: \omega_\nu^2 \leq s_l/2 \}$.

Then, with probability at least $1-\sum_{l=1}^L \gamma_l$, this approximation fulfills the bounds
\begin{align}
\label{eq:linfBound}
\|F-F^L_{\text{MLCS}}\|_\infty 
&\leq C h_L^{t+t'} \|f\|_{\Ycal_t'}\|\Gcal\|_{\Xcal_{t'}'}(L+1), 
\\
\|F-F^L_{\text{MLCS}}\|_2 
&\leq C' 
      h_L^{t+t'} \|f\|_{\Ycal_t'}\|\Gcal\|_{\Xcal_{t'}'}
\end{align}
and can be computed in a total work that scales as 
\begin{equation}
W_L^T \lesssim \left\{ \begin{array}{cc} 2^{nL}, & \sigma_p(\tau) < n, \\
L^{\xi+1} 2^{nL}, &\sigma_p(\tau) = n, \\
L^\xi 2^{L\sigma_p(\tau)}, &\sigma_p(\tau) > n, 
 \end{array} \right.
\end{equation}
where $\xi = 4$ or $5$ for constant or polynomially growing weights $\vbf$, respectively. 
\end{atheorem}

\begin{proof}
This theorem follows from applying Theorem~\ref{thm:wl1Truncated} at each level $l \in \{1,\hdots,L\}$
with probability \changes{of failure $\gamma_l$ and taking the union bound.} 
The bound~\eqref{eq:linfBound} follows from 
Theorem~\ref{thm:summability} and using the \changes{calculations in Prop.\ref{prop:computations}.}
\end{proof}

\changes{Note that we can make the failure probability more explicit. Choosing
$s_l$ as in~\eqref{eq:sl}, i.e., $s_l = \Ccal 2^{\sigma_p(L-l)}$ with
$\sigma_p = \frac{(t+t')p}{1-p}$ and equating both 
terms in the $\max$ defining the number of samples so that 
\[
m_l \asymp s_l \log(s_l)^3 \log(N_l) \asymp  \sigma_p^3 2^{\sigma_p (L-l)} (L-l) \log(N_l), \quad l=1,\hdots,L
\]
gives
\begin{equation*}
\gamma_l = N_l^{-\log^3(s_l)} = N_l^{-\(c_1 \sigma_p(L-l) + c_2\)^3}, \quad l=1,\hdots,L.
\end{equation*}
This results in a probability of failure at most
\begin{align*}
\sum_{l=1}^L \gamma_l & = \sum_{l=1}^L N_l^{-\log^3(s_l)} = \sum_{l=1}^L N_l^{-(\(c_1 \sigma_p(L-l) + c_2\)^3}
\leq \sum_{l=1}^L N_L^{-(\(c_1 \sigma_p(L-l) + c_2\)^3} \lesssim N_L^{-c_2^3} \sum_{l=1}^L N_L^{-(c_1 \sigma_p^3 (L-l))^3}\\
& \lesssim N_L^{-c_3},
\end{align*}
where it is used that $N_1 \geq N_2 \geq \cdots \geq N_L$ by definition of $\Gamma_l$ and $s_l$.}

\section{Numerical results}
\label{sec:numResults}
In this section, we illustrate our theoretical findings with some numerical examples. 
\changes{All these examples are implemented in Python 2.7. For the PDE solves we use the tools developed via the FEniCS project}~\cite{AlnaesBlechta2015a,LoggMardalEtAl2012a}.\footnote{Note that all the code for reproducible research and further use is available from one of the authors' github page: \url{https://github.com/jlbouchot/CSPDEs}}.
\changes{The sparse recovery problem is either solved by greedy approaches} using functions developed in house or using CVXPY~\cite{cvxpy} as convex solver for the weighted $\ell_1$ problem. 
\changes{We use degree one Lagrange polynomials as finite elements.} 

\subsection{Convergence}
\label{ssec:convergence}
We start by looking at the convergence of the approach with respect to the meshwidth. 
\changes{To this end}, we fix the number of levels used for the MLCSPG approximation to $L = 3$ and let the coarser meshwidth vary from $h_0 = 1/5$ to $h_0 = 1/70$. 
We want to illustrate the results suggested by Theorem~\ref{thm:mainThm}. 
\changes{We consider the diffusion problem~\eqref{eq:Diffusion} and the QoI $F(\bsy) = \int \limits_{x \in D}u(x,\bsy)\dif{x}$.  The diffusion coefficient is represented via a cosine expansion, i.e., we set} 
\begin{equation}
\label{eq:cosineExpansion}
a(x,\bsy) = \bar{a}(x) + \sum_{j=1}^d y_j \frac{\cos(\pi j \|x\|)}{j^\mu}, \changes{\quad x \in \Rbb^n, \bsy \in [-1,1]^d}.
\end{equation}
For the results presented here, we 
choose $n =1,2,3$, $\mu = 2$, and $d = 10,15,20$. 
We set uniform weights $v_j = 1.08$ and $\omega_\nu = \theta^{\|\nu\|_0}\vbf^\nu$ as suggested in Theorem~\ref{thm:mainThm}. 
The mean field and forcing term are kept constant $\bar{a} \equiv 4.3$ and $f \equiv 10$, respectively. 

\begin{figure}[htbp]
\centering
\includegraphics[width=0.85\linewidth]{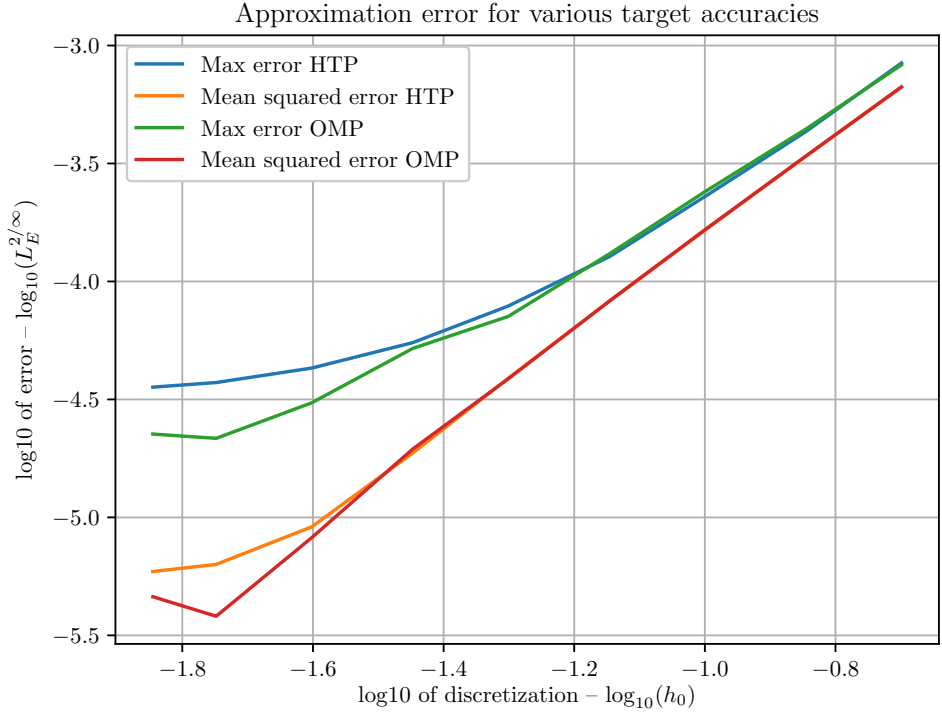}
\caption{
Convergence of the MLCSPG method for the diffusion problem with cosine expansion of the diffusion coefficient. $d=20$. The graph shows $\log(\text{error})$ as a function of $\log(h_0)$ where the empirical error is measured in the $L_E^2$ and $L_E^\infty$ norm and for sparse recovery done with either HTP or OMP. See text for more details.
 }
\label{fig:loglogplot}
\end{figure}

Figure~\ref{fig:loglogplot} illustrates the convergence of the method with respect to the coarser meshwidth $h_0$ in $n=2$ spatial dimensions and $d=20$ parameters. 
The empirical errors are calculated as 
\begin{align*}
L_E^1 &:= \frac{1}{N_{\text{test}}}
\sum_{1 \leq j \leq N_{\text{test}}} |F(\bsy^{(j)}) - F^{\text{CSPG}}(\bsy^{(j)})| \;, 
\\ 
L_E^2 &:= \sqrt{\frac{1}{N_{\text{test}}}
\sum_{1 \leq j \leq N_{\text{test}}} |F(\bsy^{(j)}) - F^{\text{CSPG}}(\bsy^{(j)})|^2} \;, 
\\
L_E^\infty &:= \max_{1 \leq j \leq N_{\text{test}}} |F(\bsy^{(j)}) - F^{\text{CSPG}}(\bsy^{(j)})| 
\;,
\end{align*}
for $N_{\text{test}} = 1000$ independent draws of random parameter vectors $\bsy$ and where the ground truth used for comparison is a numerical approximation computed on a grid that is at least $4$ times finer. 
The sparse recovery methods used for this figure are iterative (HTP ~\cite{Foucart11HTP}, picked for its proven fast convergence~\cite{Bouchot13GHTP}) and greedy (OMP) approaches. 
The finite element method is used with a degree one polynomial and with an iterative Krylov solver for the inversion of the system involving the stiffness matrix. 

The (level dependent) number of samples and sparsities have been chosen as 
\begin{align}
m_l &= 2 \cdot s_l \log(N_l), \label{eq:sprse1d_m}\\ 
s_l &= 8\cdot 2^{L-l}\;, \label{eq:sprse1d}
\end{align}
where $N_l = |\Gamma_l|$ is the set of \Tscheb polynomials truncated according to Eq.~\eqref{eq:gamma0}.
The choice of $m_l$ differs slightly from 
the theoretically justified choice in Eq.~\eqref{eq:ml}.
The selection \eqref{eq:sprse1d_m} refers to the usual rule of thumb in compressed sensing 
which is justified by non-uniform recovery results with random matrices,
see~\cite[Ch.9.2]{Foucart13book} for details. 
While the choice \eqref{eq:sprse1d_m} of numbers of CS sample $m_l$ 
is below what is sufficient according to our theoretical results, we 
shall see in the numerical examples ahead that even this optimistic selection of
sample number is more than sufficient for our problems. 
The choice \eqref{eq:sprse1d} of $s_l$ corresponds to Eq.~\eqref{eq:sl} where the 
proportionality constant is chosen as $8$ and the regularity assumption 
of the solution is taken as $\sigma_p(t+t') = (t+t')p/(1-p) = 1$ to simplify the exposition. 
This constant could be estimated numerically from Figure~\ref{fig:loglogplot}. 

\subsection{Computing times}
\label{ssec:compTime}
We now investigate the actual computational complexity required for our approach. 
We consider the same framework as in the previous section but consider $d = 15$ parameters in $n = 2$ and $3$ spatial dimensions. 
All other parameters are kept the same. 
\begin{figure}[htbp]
\centering
\includegraphics[width=0.85\linewidth]{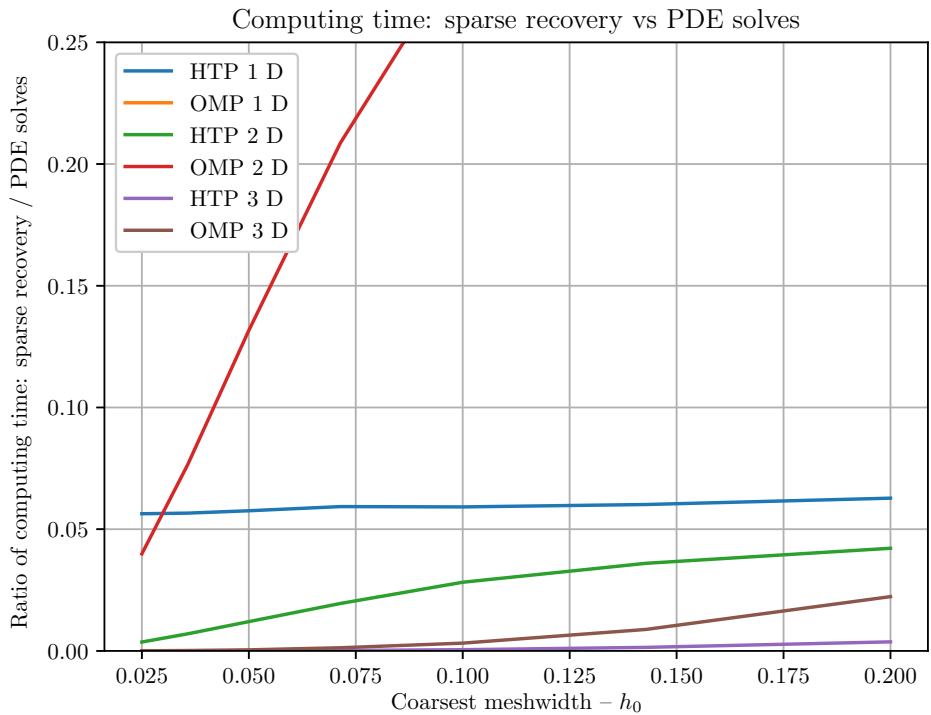}
\caption{Ratio of time for the computations for the sparse recovery against the time for the PDE solves. 
As the model gets more and more complex, the computing time for the sparse recovery is more and more negligible.  }
\label{fig:sparseagainstsolves}
\end{figure}

Figure~\ref{fig:sparseagainstsolves} shows the ratio of the time required for the sparse recovery problems against the time required for computing the PDE solutions. 
The times (for sparse recovery and PDE solves) are reported by adding all contributions at every level and then the ratio sparse recovery to PDE solves is graphed. 
We displayed the results when using weighted versions of OMP and HTP which should be prefered in high-dimensions. 
As the complexity increases (red and green curves for the case of $n=3$ spatial dimensions) the time required for the sample evaluations becomes more and more important compared to the time needed for the sparse recovery. 
Note that the orange curve displaying the computing time for the weighted OMP in $n = 2$ dimensions, while being cut at the top, never reaches more than $0.35$ in our experiments.
Table~\ref{tab:computingTime} shows a precise description of the time required for the recovery and sample evaluations at each level when considering $d = 15$ parameters and two different original meshsizes $h_0 = 1/5$ and $h_0 = 1/40$. 
As it can be seen, the efficiency of the presently proposed MLCSPG approach increases for more expensive forward solves, i.e. with increasing complexity of the simulated system. 
Moreover, it is important to point out that the sparse recovery component is completely independent of the size of the spatial dimension and its discretization as illustrated in Table~\ref{tab:computingTime}.

\begin{table}[]
\centering
\caption{Comparison of the different computing times for various settings (in seconds)}
\label{tab:computingTime}
\begin{tabular}{|l|l|l|l|l|l|l|l|l|l|l|l|l|}
\hline
\multirow{2}{*}{} & \multicolumn{4}{l|}{$l = 0$}                         & \multicolumn{4}{l|}{$l = 1$}  & \multicolumn{4}{l|}{$l = 2$}                         \\ \cline{2-13} 
                  & \multicolumn{2}{l|}{HTP} & \multicolumn{2}{l|}{PDE} & \multicolumn{2}{l|}{HTP} & \multicolumn{2}{l|}{PDE} & \multicolumn{2}{l|}{HTP} & \multicolumn{2}{l|}{PDE} \\ \cline{1-13} 
   $h_0 =$ & $1/5$ & $1/40$ & $1/5$ & $1/40$ & $1/5$ & $1/40$ & $1/5$ & $1/40$ & $1/5$ & $1/40$ & $1/5$ & $1/40$ \\ \hline
   $n = 1$ & 2.8876 & 2.8879 & 24.324 & 24.787 & 0.1730 & 0.1156 & 19.651 & 20.186 & 0.0111 & 0.0108 & 7.6102 & 8.0154 \\ \hline
   $n = 2$ & 2.9699 & 2.8923 & 32.247 & 158.59 & 0.2051 & 0.2048 & 27.987 & 286.88 & 0.0129 & 0.0926 & 15.339 & 428.39 \\ \hline
   $n = 3$ & 2.9079 & 2.8167 & 84.965 & 24221 & 0.2034 & 0.1946 & 207.31 & 88919 & 0.0119 & 0.0961 & 569.01 & 286804 \\ \hline
\end{tabular}
\end{table}


\subsection{Single-level versus Multi-level CSPG}
\label{ssec:SL_vs_ML}
\begin{figure}[htbp]
\centering
\includegraphics[width=0.85\linewidth]{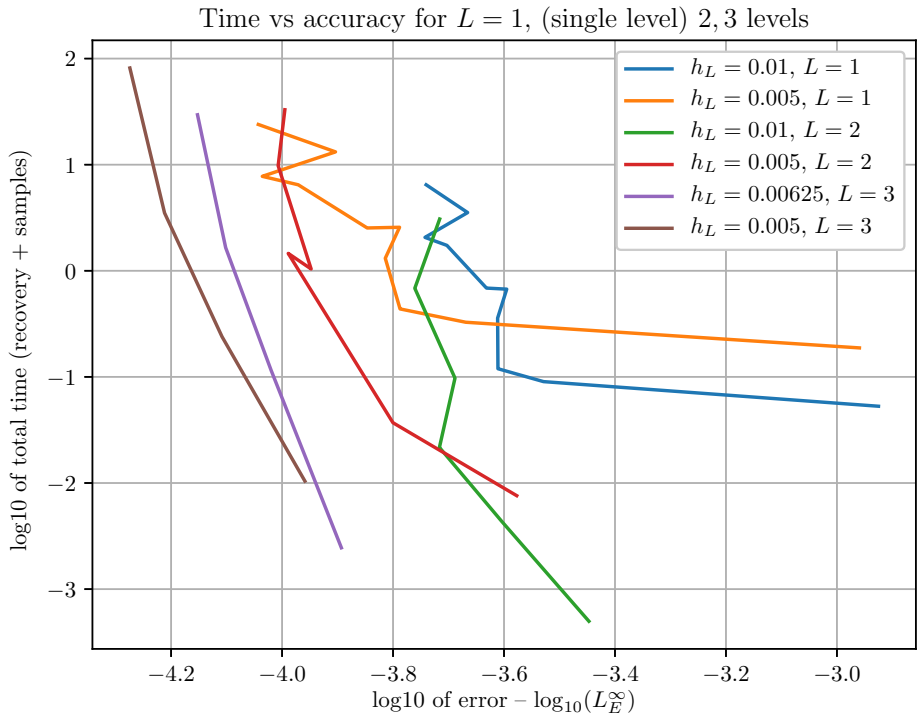}
\caption{Time for evaluating the samples and sparse recovery against accuracy for $L=1,2,3$ levels of approximations. The case $L = 1$ corresponds to the SLCSPG approach. }
\label{fig:SL_vs_ML}
\end{figure}

To ensure the necessity of the multi-level approach developed  in this paper, we compare with the original single-level method introduced by two of the named authors~\cite{Rauhut14CSPG}. 
We report on the computational time required to reach a given accuracy, parametrized by fixing $h_L$. 
We compare the single-level approach ($L=1$) with multi-level schemes $L=2,3$. 
By varying the constants of proportionality related to the sparsity per level (see Eq.~\eqref{eq:sprse1d}) we can control the computational time. 
Fig~\ref{fig:SL_vs_ML} illustrates this behavior when dealing with $d = 15$ parameters by plotting the logarithm of the computing time (including both the sample evaluations and the sparse recovery parts) with respect to the (log of the) accuracy. 

Once again we consider the parametric diffusion problem from Eq.~\eqref{eq:Diffusion} in $n = 2$ spatial dimensions and with the cosine expansion of the diffusion coefficient as in Eq.~\eqref{eq:cosineExpansion}. 
The cosine coefficients have decay parameter $\mu = 2$ and the mean-field is set to be a constant coefficient $\bar{a}(x) =  4.3$, for all $x \in [0,1]^2$. 

The proportionality constant in \eqref{eq:sprse1d} (where it is $8$) varies from $5$ to $14$ for the single-level case, from $3$ to $8$ for the $L=2$ case, and from $3$ to $6$ for the $L=3$ case. 
First, as suggested by the SLCSPG theory (see \cite{Rauhut14CSPG,Bouchot15SLCSPG}), the single level approach is indeed limited by the quality of the approximation grid and eventually reaches a maximum accuracy. 
At this point, adding more samples does not yield better approximation capabilities (see the blue and orange curves). 
On the other hand, using the multi-level approach, we can reach better accuracy while reducing the computations, illustrated by the fact that the green, red, purple and brown curves are all below and to the left of the SLCSPG curves (blue and orange curves in Fig.~\ref{fig:SL_vs_ML}).

Comparing the curves for the highest resolution plotted ($h_L = 0.005$, orange for the single level computations, red for two levels, and brown for the three level case), the advantage in terms of computations and accuracy of the presented methods becomes noticeable. 
Interestingly, the use of a multi-level procedure seems to break (to some extent) the mesh-size barrier imposed by the single-level -- a remark worth investigating further.

\subsection{Variations with the dimension}

To ensure the numerical scalability of the approach as the parametric dimension increases, we ran tests on a piecewise constant diffusion problem in two spatial dimensions. 
In this set of experiments we partition the spatial domain $[0,1]^2$ into $d = 3^2, 4^2, 5^2,$ and $6^2$ equal patches $D_j$, $1 \leq j \leq d$. 
The diffusion coefficient from Eq.~\eqref{eq:Diffusion} is given by 
\begin{equation}
a(x,\bsy) = \bar{a}(x) + v\sum_{j = 1}^d y_j\chi_{D_j}(x),
\end{equation}
where $\chi_{D_j}$ corresponds to the characteristic function of the patch $D_j$. 
We report our experiments in Table~\ref{table:2DPWCD}, where the local variations have an amplitude up to $v = 2$ and the mean field is constant $\bar{a} \equiv 5$. 
This table shows the time required for computing the PDE solutions as well as the sparse recovery procedures (done via the weighted version of HTP here) for the different numbers of parameters. 
We also mention the number of samples and the size of the truncated active set $\Gamma_l$ at every level.

\begin{table}[]
\centering
\caption{Statistics on the computing time, number of samples and size of active sets, and accuracy for a piecewise constant diffusion problem in two spatial dimensions with various parametric dimensions. These values were obtained using FEniCS as direct forward solver, and with a Python implementation. }
\label{table:2DPWCD}
\begin{tabular}{|l|l||c|c|c|c|}
\hline
\multicolumn{2}{|l||}{Size} & $3 \times 3$ & $4 \times 4$  & $5 \times 5$ & $6 \times 6$ \\ \hline
\multicolumn{2}{|l||}{$d $} & $9$ & $16$ & $25$ & $36$ \\ \hline \hline
\multicolumn{2}{|l||}{$L_E^2$} & $0.000271$ & $0.000478$ & $0.000845$ & $0.001005$ \\ \hline
\multicolumn{2}{|l||}{$L_E^\infty$} & $0.001176$ & $0.001687$ & $0.004229$ & $0.004081$ \\ \hline
\multicolumn{2}{|l||}{PDE solves (s)} & $2991.9141$ & $4304.8436$ & $6049.8758$ & $23437.15047$ \\ \hline
\multicolumn{2}{|l||}{Recovery (s)} & $0.3874$ & $1.7044$ & $5.9584$ & $29.7932$ \\ \hline \hline
\multicolumn{2}{|l||}{$\log($time $d_i)$/$\log(d_{i})$} & $3.6427$ & $3.0180$ & $2.7055$ & $2.8082$ \\ \hline \hline
\multirow{3}{*}{$m$}  & $l=0$  & $1082$ & $1298$ & $1467$ & $1605$ \\ \cline{2-6} 
                   & $l=1$  & $438$ & $518$ & $582$ & $637$ \\ \cline{2-6} 
                   &  $l=2$ & $165$ & $198$ & $224$ & $246$ \\ \hline \hline
\multirow{3}{*}{$N$}  &  $l=0$ & $4687$ & $25225$ & $94351$ & $278755$ \\ \cline{2-6} 
                   & $l=1$  & $931$ & $3241$ & $8851$ & $20731$ \\ \cline{2-6} 
                   &  $l=2$ & $172$ & $473$ & $1076$ & $2143$ \\ \hline
\end{tabular}
\end{table}

As claimed in this article, we are capable of breaking the curse of dimensionality. 
Indeed, computing the ratio of the logarithms of the computing times to the dimension shows that complexity only scales polynomially in the number of parameters (as claimed in Theorem~\ref{thm:mainThm}).
Moreover, having set a target accuracy of $h_L = 2.5 \cdot10^{-3}$ we verified the accuracy of our recovered solutions against 1000 random independent tests. 
The ground truth was here numerically approximated on a grid three times finer than the one used for the MLCSPG method. 

\subsection{Comparison with $L^2$ and Monte-Carlo}

We consider the bounded interval $D = (0,1)$ with equispaced
partition $\overline{D} = \bigcup_{i=1}^d \overline{D_i}$ into subintervals
$D_i = (x_{i-1}, x_i)$ where $x_i = i/d$, for some $d\in \IN$.
We let $a(x,\bsy) = \bar{a} + \sum_{j=1}^d y_j c_j \chi_{D_j}(x)$ 
with $\bar{a}$ being a constant independent of $x$, $\{c_j\}_{j=1}^d$ a 
predefined (fixed) sequence such that the (weighted) uniform ellipticity 
assumption~\eqref{UEA} holds, and $\chi_{D_j}$ the indicator function of the set $D_j$. 
We select the parameters to be 
$d = 6$ and pick uniform (small) local variations as $c_j = 1/6$, for $1 \leq j \leq d$.
The uniform weights $v_j$ are selected as $v_j = 1.07$ for all $j$. 
We also set the forcing term $f \equiv 1$ to be constant.
Then, \emph{for any $\bsy \in U$ the solution to the diffusion equation 
is continuous and piecewise quadratic.
}
The (level dependent) number of samples and sparsities are chosen as 
\begin{align*}
m_l &= 2 \cdot s_l \log(N_l), \\ 
s_l &= 20\cdot 2^{L-l}\;. 
\end{align*}

The initial mesh size is set to $h_0 = 5 \cdot 10^{-4}$. Further numerical tests -- not included in this paper, 
but available online -- have shown that this parameter has, in this case, 
little to no influcence over the results.

\begin{figure}[htbp]
\centering
\includegraphics[width=0.75\linewidth]{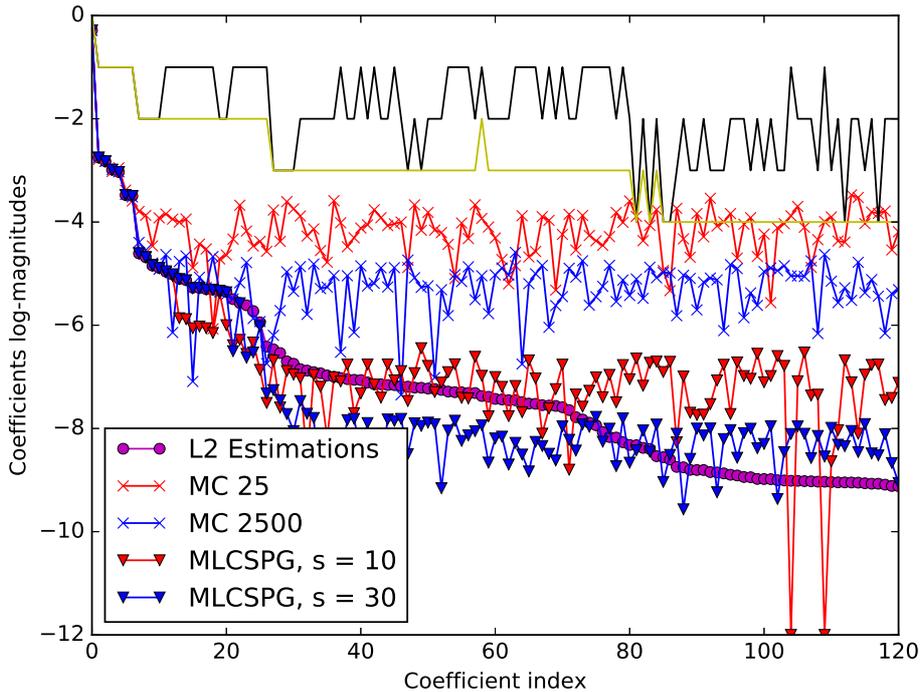}
\caption{
Estimations of the 120 largest coefficients (in magnitude) of the gpc of the piecewise constant diffusion problem (see text for details on the parameters)
reordered by decreasing magnitude, according to their estimations via a least squares method. 
The {\color{purple}{least squares estimation}} computed $12171$ coefficients from $24342$ 
random samples; 
the {\color{blue}{downward blue triangles}} correspond to a constant $15$ in Eq.~\eqref{eq:sprse1d},
while the {\color{red}{downward red triangles}} correspond to the constant $5$. 
The {\color{red}{$MC~25$ curve}} corresponds to Monte Carlo estimations with $2.5\cdot10^{7}$
samples 
while the {\color{blue}{$MC~2500$ curve}} is based on $2.5\cdot10^{9}$ samples. }
\label{fig:coefEstimations}
\end{figure}

We compare the convergence of our algorithms with other methods: 
Monte-Carlo sampling and least squares ($\ell_2$ recovery) \cite{Migliorati14discreteL2}. 
The estimation of the \Tscheb coefficients are displayed in Fig.~\ref{fig:coefEstimations}, 
where the magnitudes of the Chebyshev coefficients of the (functional of the) 
parametric solution are displayed on a logarithmic (base 10) scale. 
The $x$-axes corresponds to an enumeration of the multi-index of the Chebyshev coefficient, 
whereby the larger ones (in magnitude, according to the $\ell_2$ recovery) are first.
The least squares solution is obtained as follows. 
We first build the active set of candidates for the truncated polynomial space 
as predicted by Theorem \ref{thm:wl1Truncated}, i.e.
$\Gamma = \cup_{l=1}^3 \Gamma_l$. 
This set has total dimension $N = 12171$. 
Then $m = 24342$  sampling points $\bsy^{(i)}$ are chosen at random, and the values $b_i = F(\bsy^{(i)})$ are computed and stacked into a vector $\bsb = (b_i)_i$. 
Finally, the coefficients $(F^{\ell_2}_\nu)_{\nu \in \Gamma}$ are computed 
as the minimizer of the least squares problem
\[
\min_F \|\bsb - \Abf F\|_2,
\]
where $\Abf_{i,\nu} = T_\nu(\bsy^{(i)})$, with $1 \leq i \leq m$ and $\nu \in \Gamma$.
To display our results on Fig.~\ref{fig:coefEstimations}, we have an (implicit) enumeration $\pi: \{1,\hdots, 12171\} \to \Gamma$ such that $|F^{\ell_2}_{\pi(1)}| \geq |F^{\ell_2}_{\pi(2)}| \geq \cdots \geq |F^{\ell_2}_{\pi(12171)}|$. 
For this experiment, we compute the solutions to the weighted 
$\ell_1$-minimization problems using the SCP optimization procedure from the 
CVXPY package~\cite{cvxpy} with accuracy for the numerical optimization  set to $10^{-6}$.
The downward triangles are the results using our suggested 
MLCSPG method with the multiplicative constant $8$ in Eq.~\eqref{eq:sprse1d} replaced by $5$ (red curve) and by $15$ (blue curve).
The selection of the constant equal to $15$ 
corresponds to $m_1 = 2258$ solves at the coarsest level, 
$m_2= 968$ at the second, and to $m_3 = 394$ solves at the finest discretization level $L=3$.
$m_1 = 576$, $m_2 = 228$, and $m_3 = 73$ samples, for the red curve. 
The crosses correspond to the MC simulations, where we have used $m = 2.5 \cdot 10^{7}$ (red curve) 
and $m = 2.5 \cdot 10^{9}$ (blue curve) samples for the estimation of the \Tscheb coefficients. 
Noting that the $y$ values of the graphs correspond to the $\log_{10}$ 
of the magnitude of the coefficients, we see that the accuracy of the MC estimations is limited by the
mean square convergence rate $m^{-1/2}$.
The purple circles correspond to the $\ell_2$ estimation described above
(this corresponds to an oversampling ratio of $2$, which is far below theoretical results).
In this example, our approach (as illustrated by the downward triangle curves in Fig.~\ref{fig:coefEstimations}) produces reliable approximations of gpc coefficients which are large in magnitude 
with a number of samples orders of magnitudes smaller than both the $\ell_2$ and the MC approaches.
The limitation of the MC method to a square root convergence rate 
requires a prohibitive number of samples for more complicated PDEs.
It is important to note also that the accuracy of the recovered coefficients 
via our MLCSPG method are constrained by the accuracy of the numerical solver for the weighted $\ell_1$ minimization. 
Finally, the yellow curve corresponds to the (negative, for illustrative purposes) total degree of the multi-index of the associated \Tscheb coefficient while the black curve corresponds to the (negative of the) maximum degree in the tensor product~\eqref{eq:TensTscheb}. 
It is interesting to notice that the magnitude seems to be smaller as the degree of the multi-index increases.

\section{Conclusions}
\label{sec:Concl}
For a class of abstract, 
affine-parametric, linear operator equations depending on sequences $\bsy$
of parameters, we have introduced a multi-level generalization of the 
CS approach from \cite{Rauhut14CSPG} to efficiently 
scan the high-dimensional parameter space. 
For the approximate solution of (instances of) the parametric
operator equations, we stipulated available inf-sup stable, 
Petrov-Galerkin (``PG'' for short) discretizations of the 
``nominal'' operator $A_0 = A(\bsnul)$; 
in particular, \eqref{eq:DiscInfSupNom} holds.
The small perturbation hypothesis \eqref{eq:sumbjt} at $t=0$
implies uniform (w.r.t. $\bsy\in U$) inf-sup stability
\eqref{eq:infsupA} of the PG discretization (Thm. \ref{thm:uniinfsup}). 
Admissible PG discretizations comprise, in particular, 
all classical primal or mixed Finite Element Methods
(FEM for short), as well as spectral and collocation methods 
for elliptic and certain linear, parabolic evolution equations.
Throughout, we used multi-level
Finite Element Galerkin discretizations in $D\subset \IR^n$ with 
isotropic mesh refinements, 
responsible for the $\mathcal{O}(2^{nl})$ scaling in the proof of 
Proposition \ref{prop:computations}. 
Anisotropic, ``sparse-grid'' discretizations of the parametric problems in $D$ 
would result, with analogous analysis, in
so-called ``multi-index'' compressed sensing PG methods,
analogous to multi-index MC in \cite{HAFNRT_MIMC2016}, with
$\mathcal{O}(l^{n-1}2^{l})$ in place of $\mathcal{O}(2^{nl})$.
We analyzed error vs. work of the multi-level extension of the 
combined, CS-PG algorithm and showed that it affords improved,
as compared to the single-level variant from \cite{Rauhut14CSPG,Bouchot15SLCSPG}, 
error vs. work bounds with convergence rates that are independent
of the dimension of the space parameters which are active in the
approximation, while being \emph{``nonintrusive''}, i.e.
accessing an available solver at each discretization level.
This is analogous to what is known from multi-level Monte-Carlo 
(``MLMC'' for short) sampling methods, 
as surveyed e.g. in \cite{MLMCGilesActa}.
Contrary to MLMC methods whose convergence rate is limited
by the (mean-square) rate $1/2$ afforded by MC methods, 
and the recently proposed sparse-grid methods 
in \cite{chkifa2014polyDownward} which rely on a particular 
(``downward closed'') structure of the sets of active polynomials, however,
the presently proposed approach yields dimension-independent
convergence rates (potentially far beyond $1/2$) in the sup-norm with respect to
the parameters, exploiting any sparsity 
in the gpc coefficient sequence of the parametric solutions,
\emph{without} strong, a-priori structural 
assumptions on the active polynomial degrees.
At the same time, 
the MLCSPG approach is nonintrusive and intrinsically parallel 
as MLMC methods. 
If a-priori information on the structure of sets of active 
indices (such as ``downward closedness'') is available, 
corresponding accelerations of the SLCS approach have
recently been investigated in \cite{ChkifDextTranWebs16}.
This is afforded by adopting
\Tscheb gpc expansions which are orthonormal with respect to a
probability measure which underlies the CS method, whereas
sparse-grid methods as in \cite{chkifa2014polyDownward}
afford greater flexibiliy as regards the choice of gpc system.

We remark that although here only affine-parametric operator equations
were considered, the key results of the present paper require merely
\emph{sparsity of \Tscheb gpc expansions} 
(as expressed, e.g., in summability
of sequences of $\Xcal_t$-norms of gpc expansion coefficients 
in the conditions \eqref{eq:sumbjvjtilda} - \eqref{eq:admissibility},
rather than the weaker summability of $\Xcal$-norms 
in the SLCSPG considered in \cite{Rauhut14CSPG}) 
of the parametric solutions, and some (possibly crude) bounds of these
coefficients which enter the weight sequence $\omega$, 
and a family of uniformly inf-sup stable PG discretization methods.
Such results are available for rather general,
holomorphic-parametric, nonlinear operator equations 
in \cite{Chkifa14Breakingcurse}.
In case that the $\psi_j$ in \eqref{eq:aKL}, \eqref{eq:Diffusion} 
have supports which are localized to subdomains of $D$ with controlled
overlap, higher summability for the \Tscheb gpc expansion coefficients
holds; we refer to \cite{Bachmayr2015sparsepoly} for details.
The presently proposed MLCSPG algorithms are able to exploit 
better summability of \Tscheb gpc expansion coefficients without any modification
in the algorithm.

\bibliographystyle{abbrv}
\bibliography{MLCSPG}
\end{document}